\documentclass[12pt]{elsarticle}
\usepackage{amsthm, amsmath,amssymb,amsfonts}
\usepackage{graphicx}
\usepackage{color}
\usepackage{psfig}
\usepackage{graphics}
\usepackage{amsmath}
\usepackage{amssymb}
\usepackage{epsfig}
\usepackage{epstopdf}
\usepackage{color}

 \newtheorem{proposition}{Proposition}
\newtheorem{theorem}{Theorem}
\newtheorem{corollary}{Corollary}
\newtheorem{lemma}{Lemma}
\newtheorem{remark}{Remark}
\newtheorem{definition}{Definition}
\newtheorem{example}{Example}

\def\para{\vspace{1.5mm}}
\def\cP{{\mathcal P}}

\def\gcd{{\rm gcd}}

\def\lc{{\rm lc}}
\def\Content{{\rm Content}}

\def\mult{{\rm mult}}

\def\deg{{\rm deg}}
\def\degree{{\rm deg}}

\def\Res{{\rm Res}}

\begin{document}

\begin{frontmatter}

\title{Resultants and Singularities of Parametric Curves}
\author{Angel Blasco and Sonia P\'erez--D\'{\i}az\\
Departamento de F\'{\i}sica y Matem\'aticas \\
        Universidad de Alcal\'a \\
      28871-Alcal\'a de Henares, Madrid, Spain  \\
angel.blasco@uah.es, sonia.perez@uah.es
}
%\date{}          % Enter your date or \today between curly braces
%\maketitle

\begin{abstract} Let ${\cal C}$ be an algebraic space curve defined parametrically by $\cP(t)\in {\Bbb K}(t)^{n},\,n\geq 2$. In this paper, we introduce a polynomial,  the {\it T--function}, $T(s)$, which is defined by means of a univariate resultant constructed from $\cP(t)$. We show that $T(s)=\prod_{i=1}^n H_{P_i}(s)^{m_i-1}$, where $H_{P_i}(s),\,i=1,\ldots,n$ are polynomials (called {\it the fibre functions}) whose roots are the fibre of the ordinary singularities  $P_i\in {\cal C}$ of multiplicity $m_i,\,i=1,\ldots,n$. Thus, a complete classification of the singularities of a given space curve, via the factorization of a resultant, is obtained.
\end{abstract}


\end{frontmatter}

\section{Introduction}

Parametrizations of rational curves play an important role in many practical applications in computer
aided geometric design where objects are often given and manipulated
parametrically  (see e.g. \cite{Hoffmann}, \cite{HSW}, \cite{HL97}). In the last years, important advances have been made concerning the information one may obtain from a given rational parametrization defining an algebraic variety. For instance, a complete analysis of the asymptotic behavior of a given curve has been carried out in \cite{Blas2}; efficient algorithms for computing the implicit equations that define the curve are provided in  \cite{Buse2010} and \cite{SWP} and the study and computation of the fibre of a point via the parametrization  can be found in  \cite{SWP}. In addition, some aspects concerning the singularities of the curve and their multiplicities are studied in \cite{Abhy}, \cite{Buse2012}, \cite{Chen2008}, \cite{JSC-Perez} and \cite{Rubio}. Similar problems, for the case of a given rational parametric surface, are being analyzed. For instance, the computation of the singularities and their multiplicities from the input parametrization is presented in \cite{MACOM}, a univariate resultant-based implicitization algorithm for surfaces is provided in \cite{Impli-Super}, and the computation of the fibre of rational surface parametrizations is developed in \cite{Fibra-Super}.

\para

In this paper, we show how to relate the fibre  and the singularities of a given curve defined parametrically, by means of a univariate resultant which is constructed directly from the parametrization. For this purpose, we consider
$\cP(t)\in\mathbb{P}^n(\mathbb{K}(t))$ a rational projective parametrization of an algebraic curve $\cal C$ over an algebraically closed field of characteristic zero, $\Bbb K$.  Associated with   ${\mathcal P}(t)$, we
consider the induced rational map $\psi_{\mathcal
P}:{\mathbb{K}}\longrightarrow {\mathcal C} \subset \mathbb{P}^n(\mathbb{K}); t\longmapsto
{\mathcal P}(t).$ We denote by $\degree(\psi_{\mathcal P})$ the degree of the rational map
$\psi_{\mathcal P}$. The birationality of  $\psi_\cP$, i.e. the properness of $\cP(t)$, is
characterized  by $\deg(\psi_\cP)=1$ (see \cite{Harris:algebraic}
and \cite{shafa}). Intuitively speaking, $\cP(t)$ proper means that $\cP(t)$ traces the curve once, except for at most a finite number of points. We will see that, in fact,  these points are the singularities of $\cal C$.

We recall that the degree of a rational
map can be seen as the cardinality of the fibre of a generic element
(see \cite{shafa}). We  use this
characterization in our reasoning and thus, we denote by
${\mathcal F}_{\mathcal P}(P)$ the fibre of a point $P\in \mathcal
C$ via the parametrization $\mathcal{P}(t)$; that is $ {\mathcal F}_{\mathcal P}(P)={\mathcal P}^{-1}(P)=\{
t\in \mathbb{K} \,|\, {\mathcal P}(t)=P \}. $

\para

In order to make   the paper more reader--friendly, we first consider the case of a given plane curve $\cal C$ defined parametrically by  $\cP(t)\in \mathbb{P}^2(\mathbb{K}(t))$ (see Sections  \ref{S-notacion} and \ref{S-resultantsingularities})  to, afterwards, generalize the results obtained to rational space curves in any dimension (see Section \ref{S-spacialcurves}). We also assume that $\cal C$ has only ordinary singularities (otherwise, one may apply quadratic transformations for birationally transforming the curve into a curve with only ordinary singularities). Non--ordinary singularities have to be treated specially since a non--ordinary singularity might have other
singularities in its ``neighborhood''. This specific case will be addressed in a future work and in fact, we will show that   similar results to those presented in this paper  can be stated   for curves with non--ordinary singularities.

\para

Under these conditions, the main goal of the paper is to prove that a univariate resultant constructed directly from $\cP(t)$, which we will call the  {\it T--function}, $T(t)$, describes totally the singularities of $\cal C$. It will be proved that the factorization of $T(t)$ provides the fibre functions of the different singularities of  $\cal C$ as well as their corresponding multiplicities. The {\it fibre function} of a point $P\in {\cal C}$ via $\cP(t)$ is given by a polynomial $H_P(t)$ which satisfies that  $t_0\in \mathcal{F}_\mathcal{P}(P)$ if and only if $H_P(t_0)=0$. In \cite{JSC-Perez}, it is proved that if $H_P(t)=\prod_{i=1}^n(t-s_i)^{k_i}$ then, $\mathcal{C}$ has  $n$ tangents at $P$ of multiplicities   $k_1,\ldots,k_n$, respectively. In addition, these tangents can be computed using $\cP(t)$ and the roots of each corresponding fibre function. Furthermore, it is shown that $\mult_P(\mathcal{C})=\deg(H_P(t)).$

\para

Taking into account these previous results, in this paper we prove  that the T--function can be factorized as $T(t)=\prod_{i=1}^n H_{P_i}(s)^{m_i-1}$, where $H_{P_i}(t)$ is the fibre function  of the ordinary singularity  $P_i\in {\cal C}$ and $m_i$ is its multiplicity (for $i=1,\ldots,n$). Thus, a complete classification of the singularities of a given rational curve, via the factorization of a univariate resultant, is obtained.

\para

On finishing this work, we just found a paper by Abhyankar (see \cite{Abhy}) that proves the factorization of the T--function for a given polynomial parametrization. In addition, Bus\'e et al., in \cite{Buse2012}, provide a generalization of Abhyankar's formula for the case of rational parametrizations (not necessarily polynomial). This approach is based on the concept of singular factors introduced in  \cite{Chen2008}, and it involves the construction of  $\mu$--basis. Our approach is totally different, since we generalize Abhyankar's formula by using the methods and techniques presented in \cite{JSC-Perez}. This allows us to group the factors of the T--function to easily obtain the fibre functions of the different singularities. In addition, we show how to deal with singularities that are reached by algebraic values of the parameter.

\para

As we mentioned above, these results can be stated similarly for the case of rational space curves in any dimension. We remark that the methods developed in this paper generalize some previous  results that   partially approach the computation and analysis of singularities for rational parametrized curves  (see e.g. \cite{Buse2010}, \cite{JSC-Perez} or \cite{Rubio}). Moreover, the ideas presented open several important ways that may be used to obtain significant results concerning rational parametrizations of surfaces.  In a future work, this problem will be developed in more detail and some important results are expected to be provided.

%\para

%\textcolor{red}{We observe that these results can not be applied to the  {\it limit point}, $P_L$, of the given parametrization $\mathcal{P}(t)$. This point has to be treated specially since it is not correctly represented by %$\mathcal{P}(t)$ and  as a consequence, its multiplicity is not given by the cardinality of its fibre. In \cite{MyB-2017(b)}, we obtain some remarkable properties concerning $P_L$ and, in particular, we prove that part of its %multiplicity can not be obtained via $\cP(t)$.}

%\para

%{\color{blue}In \cite{Buse2010}, the authors use $\mu$-basis to obtain an extension of the singular factors of \cite{Chen2008} and prove that the fibre functions of the singular points divide those factors. However, some other residual elements may appear that do not correspond to any fibre function.}

\para

The structure of the paper is as follows. Sections \ref{S-notacion} and \ref{S-resultantsingularities} are devoted to the study of plane curves. In particular, in Section \ref{S-notacion}, we introduce
the  terminology that will be used throughout this paper as well as some previous results.  In Section \ref{S-resultantsingularities}, we introduce the T--function and we present the main result of the paper. It claims that the factorization of the T--function provides the fibre functions of the different singularities of the curve. The proof of this result as well as some previous technical lemmas appear in Section \ref{S-proof}. Section
\ref{S-spacialcurves} is devoted to generalize the results in Section \ref{S-resultantsingularities} to parametric space curves in any dimension.  Throughout the whole paper, we
outline all the results obtained with illustrative examples.

%We remark that from the fibre function of a point, one may determine its multiplicity as well as its fibre and the tangent lines of the curve at that point.

\section{Analysis and computation of the fibre}\label{S-notacion}

Let $\mathcal{C}$ be a rational (projective) plane curve defined by the projective parametrization
$$\mathcal{P}(t)=(p_1(t):p_2(t):p(t))\in\mathbb{P}^2(\mathbb{K}(t)),$$
where $\gcd(p_1,p_2,p)=1$, and  $\Bbb K$ is an algebraically closed field of characteristic zero . We assume that
$\mathcal{C}$ is not a line (a line does not have multiple points). Let $d_1=\deg(p_1)$,
$d_2=\deg(p_2)$, $d_3=\deg(p)$, and $d=\max\{d_1,d_2,d_3\}$. Thus, we may write  $p_1$, $p_2$ and $p$ as
$$\left\{\begin{array}{l}p_1(t)=a_0+a_1t+a_2t^2+\cdots+a_{d}t^{d}\\p_2(t)=b_0+b_1t+b_2t^2+\cdots+b_{d}t^{d}\\p(t)=c_0+c_1t+c_2t^2+\cdots+c_{d}t^{d}.\end{array}\right.$$

Associated with  ${\mathcal P}(t)$, we consider the induced rational map $\psi_{\mathcal
P}:{\mathbb{K}}\longrightarrow {\mathcal C} \subset \mathbb{P}^2(\mathbb{K}); t\longmapsto
{\mathcal P}(t).$ We denote by $\degree(\psi_{\mathcal P})$ the degree of the rational map
$\psi_{\mathcal P}$ (for further details  see e.g. \cite{shafa}
pp.143, or \cite{Harris:algebraic} pp.80). As an important result,
we recall that  the birationality of  $\psi_\cP$, i.e. the properness of $\cP(t)$, is
characterized  by $\deg(\psi_\cP)=1$ (see \cite{Harris:algebraic}
and \cite{shafa}).  Also, we recall that the degree of a rational
map can be seen as the cardinality of the fibre of a generic element
(see Theorem 7, pp. 76 in \cite{shafa}). We will use this
characterization in our reasoning. For this purpose, we denote by
${\mathcal F}_{\mathcal P}(P)$ the fibre of a point $P\in \mathcal
C$ via the parametrization $\mathcal{P}(t)$; that is $$ {\mathcal F}_{\mathcal P}(P)={\mathcal P}^{-1}(P)=\{
t\in \mathbb{K} \,|\, {\mathcal P}(t)=P \}. $$

In general, it holds that $P\in\mathcal{C}$ if and only if $\mathcal{F}_\mathcal{P}(P)\neq\emptyset$, although an exception can be found for {\it the limit point} of the parametrization.

\para

\begin{definition}\label{D-pto-limite}
We define  {\it the limit point} of the parametrization $\mathcal{P}(t)$ as
$$P_L=\lim_{t\rightarrow\infty}\mathcal{P}(t)/t^d=(a_d:b_d:c_d).$$
\end{definition}

\para

Note that $P_L\in {\cal C}$ since  $\mathcal{P}(t)/t^d=\mathcal{P}(t)\in\mathcal{C}$, for $t\in\mathbb{K}$, and $\mathcal{C}$ is a closed set. Furthermore, we observe that, given a parametrization $\cP(t)$, there always exists an associated limit point, and it is unique.% Under these conditions, we say that a given point $P\in {\cal C}$ is a {\it  normal point} if $P\not=P_L$.

\para

The limit point is reachable via the parametrization $\cP(t)$, if there exists  $t_0\in {\Bbb K}$ such that $\mathcal{P}(t_0)=P_L$. However, the value $t_0\in {\Bbb K}$  could not exist, and then  $\mathcal{F}_\mathcal{P}(P_L)=\emptyset$. Taking into account this statement, if $P_L$ is not an affine point or it is a reachable affine point, we have that $\cP(t)$ is a {\it normal parametrization}. Otherwise, we say that $\cP(t)$ is not normal and $P_L$ is the {\it critical point} (see Subsection 6.3 in \cite{SWP}). Further properties of the limit point are stated and proved in \cite{MyB-2017(b)}.

\para

In Subsection 2.2. in \cite{SWP}, it is stated that  the degree of a dominant rational map between
two varieties of the same dimension is the cardinality of the fiber of a generic element.
Therefore, in the case of the mapping $\psi_\cP$, this implies that almost all points of $\cal C$ (except at most a finite number of points)
are generated via $\cP(t)$ by the same number of parameter values, and this number is
the degree of $\psi_\cP$. Thus, intuitively speaking, the degree measures the number of times the
parametrization traces the curve when the parameter takes values in $\Bbb K$. Taking into
account this intuitive notion, the degree of
the mapping $\psi_\cP$ is also called {\it the tracing index} of $\cP(t)$. In order to compute the tracing index, the following polynomials are considered,
\begin{equation}\label{Eq-fibra-generica}\left\{\begin{array}{l}G_1(s,t):=p_1(s)p(t)-p(s)p_1(t)\\G_2(s,t):=p_2(s)p(t)-p(s)p_2(t)\\G_3(s,t):=p_1(s)p_2(t)-p_2(s)p_1(t)\end{array}\right.\end{equation}
and $G(s,t)=\gcd(G_1(s,t),G_2(s,t),G_3(s,t))$. In the   following theorem, we compute the  tracing index  of $\cP(t)$ using the polynomial $G(s,t)$ (see Subsection 4.3 in \cite{SWP}).

\para

\begin{theorem}\label{T-tracing-index}
It holds that
$\deg(\psi_\mathcal{P})=\deg_t(G).$
\end{theorem}

\begin{remark}\label{R-grado-Gs} We observe that:
\begin{enumerate}
 \item  The polynomials $G_1$, $G_2$ and $G_3$ satisfy that  $G_i(s,t)=-G_i(t,s)$. Clearly, $G(s,t)$  also has this property.
 \item  Taking into account the above statement, it holds that $\deg_s(G_i)=\deg_t(G_i)$ for $i=1,2,3$, and $\deg_s(G)=\deg_t(G)$.
\item  It holds that $\deg_t(G_1)=\max\{d_1,d_3\}$. Indeed: if $d_1\neq d_3$, the statement trivially holds. If $d_1=d_3$,    $\deg_t(G_1)$ may decrease if $p_1(s)c_d-p(s)a_d=0$. But this would imply that $\mathcal{C}$ is
 a line, which is impossible by the assumption. Similarly, it holds that  $\deg_t(G_2)=\max\{d_2,d_3\}$, and
$\deg_t(G_3)=\max\{d_1,d_2\}$.
\item It holds that $$G(s,t)=\gcd(G_1(s,t),G_2(s,t)).$$
Indeed: since $p(t)G_3(s,t)=p_2(t)G_1(s,t)-p_1(t)G_2(s,t)$, if $h(s,t)\in {\Bbb K}[s,t]$ divides to $G_1(s,t)$ and $G_2(s,t)$, then  $h(s,t)$ divides to $G_3(s,t)$ or $p(t)$. However, if $h(s,t)$  divides $p(t)$, then $h(s,t)=h(t)$ which would imply that there exists $t_0\in\mathbb{K}$ such that
$G_1(s,t_0)=G_2(s,t_0)=p(t_0)=0$. Hence, $p_i(s)/p(s)\in {\Bbb K},\,i=1,2,$ and $\cal C$ would be a line, which is impossible by the assumption. Similarly, it holds that
$$G(s,t)=\gcd(G_1(s,t),G_3(s,t))=\gcd(G_2(s,t),G_3(s,t)).$$
\end{enumerate}
\end{remark}

\para

Throughout this paper, we assume that  $\mathcal{P}(t)$ is proper, that is  $\deg(\psi_\mathcal{P})=1$. Otherwise, we can reparametrize the curve using, for instance, the results in \cite{Perez-Repara1}.  Under these conditions, it holds that the degree of $\cal C$ is $d$ (see Theorem 6 in \cite{JSC-Perez}). In addition,    $G(t,s)=t-s$ (see Theorem \ref{T-tracing-index}) and the cardinality of the fibre for a generic point of $\cal C$ is $1$, although for a particular point it can be different. %para las singulariddes es mayor, y para el punto limite menor %{\color{red}, although the cardinal of the fibre in a particular point can be greater to $1$.}

\para

In order to analyze these special points, in the following, we consider a particular point $P=(a,b,c)\in {\cal C}$. The fibre of $P$ consists of the values $t\in\mathbb{K}$ such that $\mathcal{P}(t)= P$, that is, those which satisfy the {\it fibre equations}, defined as
\begin{equation}\label{Eq-fibre-equations}\left\{\begin{array}{l} \phi_1(t):=ap(t)-cp_1(t)=0\\ \phi_2(t):=bp(t)-cp_2(t)=0\\\phi_3(t):=ap_2(t)-bp_1(t)=0.\end{array}\right.\end{equation}

%\para

%Since  $\mathcal{F}_\mathcal{P}(P)$ are the values of
%$t\in\mathbb{K}$ such that $\mathcal{P}(t)= P$ then,   $\mathcal{F}_\mathcal{P}(P)$ is given by the values of $t\in\mathbb{K}$  that satisfy the {\it fibre equations}, defined as
%\begin{equation}\label{Eq-fibre-equations}\left\{\begin{array}{l} \phi_1(t)=ap(t)-cp_1(t)=0\\ \phi_2(t)=bp(t)-cp_2(t)=0\\\phi_3(t)=ap_2(t)-bp_1(t)=0.\end{array}\right.\end{equation}

\para

Hence, the fibre of  $P$ is given by the common roots of these equations, which motivates the following definition:

\para

\begin{definition}\label{D-fibrefunction}
Given $P\in\mathbb{P}^2(\mathbb{K})$ and the rational parametrization
$\mathcal{P}(t)\in\mathbb{P}^2(\mathbb{K}(t))$, we define  {\it the fibre function} of $P$ at
$\mathcal{P}(t)$ as
$$H_P(t):=\gcd(\phi_1,\phi_2,\phi_3).$$
 Thus,  $t_0\in \mathcal{F}_\mathcal{P}(P)$ if and only if $H_P(t_0)=0$.
\end{definition}

\para

\begin{remark}\label{R-afin-noafin} Depending on whether $P$ is an affine point or an infinity point, the   fibre function can be expressed as follows:
\begin{itemize}
\item If $P$ is an affine point, then  $c\neq 0$. Thus, $\phi_3$ can be obtained from $\phi_1$  and $\phi_2$ and, therefore, $H_P(t)=\gcd(\phi_1(t),\phi_2(t)).$
\item  If $P$ is an infinity point, then $c=0$. Thus, $\phi_1$  and $\phi_2$ are equivalent to $p(t)=0$ (note that $a\not=0$ or $b\not=0$) and, therefore,  $H_P(t)=\gcd(p(t),\phi_3(t)).$
\end{itemize}
Note that the functions $\phi_1$, $\phi_2$ and $\phi_3$ depend on $P$
 and $\mathcal{P}(t)$. However, for the sake of simplicity, we do not represent this fact in the notation.
\end{remark}

\para

In the following, we show how the fibre function  of $P$ is related with the tangents of $\cal C$ at   $P$, and with the multiplicity of $P$. For this purpose, we first recall that    $P$ is a {\it point of
multiplicity $\ell$} on ${\mathcal C}$  if and only if all the
derivatives of $F$ (where $F$ denotes the implicit polynomial defining $\cal C$) up to and including those of $(\ell-1)$--th order, vanish at
$P$ but  at least one $\ell-$th derivative does not vanish at $P$.
We denote it by $\mult_{P}({\mathcal C})$. The point $P$  is called a  {\it simple  point}  on ${\mathcal C}$  if and only if
$\mult_{P}({\mathcal C})=1$. If $\mult_{P}({\mathcal C})=\ell>1$,
then we say that $P$ is a
 multiple or   singular point (or  singularity)  of  multiplicity
$\ell$ on ${\mathcal C}$  or an   $\ell$--fold point. Clearly
$P\not\in {\mathcal C}$ if and only if $\mult_{P}({\mathcal C})=0$.

\para

 Observe that the  multiplicity of ${\mathcal C}$ at
$P$ is given as the order of the Taylor expansion of $F$ at $P$. The
tangents  to ${\mathcal  C}$ at $P$  are the irreducible factors of
the first non--vanishing form in the Taylor expansion of $F$ at $P$,
and the  multiplicity  of a tangent  is the  multiplicity of the
corresponding factor. If all the $\ell$ tangents at the $\ell$-fold point $P$ are
different, then this singularity is  called {\it ordinary}, and
{\it non--ordinary} otherwise. Thus, we say that the character of $P$ is
either ordinary or non-ordinary.

\para

In  \cite{JSC-Perez}, it is shown how to compute the singularities and its corresponding multiplicities from a given parametrization defining a rational plane curve. Furthermore, it is provided a method for computing the tangents and for analyzing the
non--ordinary singularities. In particular, the following theorem and corollary are proved.% in \cite{JSC-Perez}.

\para

\begin{theorem}\label{T-tangentes}
Let $\mathcal{C}$ be a rational algebraic curve defined by a proper pa\-ra\-me\-tri\-za\-tion $\mathcal{P}(t)$, with limit point $P_L$. Let $P\neq P_L$ be a point of $\mathcal{C}$ and let $H_P(t)=\prod_{i=1}^n(t-s_i)^{k_i}$ be its fibre function (under $\mathcal{P}(t)$). Then, $\mathcal{C}$ has  $n$ tangents at $P$ of multiplicities   $k_1,\ldots,k_n$, respectively.
\end{theorem}

\begin{remark}
It can  not be ensured that two different values of $t$, namely $s_{i_0}$ and $s_{i_1}$, provide different tangents. Thus, we could have a same tangent (at $\cP(s_{i_0})=\cP(s_{i_1})$) of multiplicity  $k_{i_0}+k_{i_1}$.
\end{remark}

\para

\begin{corollary}\label{C-multiplicidad}
Let $\mathcal{C}$ be a rational algebraic curve defined by a proper pa\-ra\-me\-tri\-za\-tion $\mathcal{P}(t)$, with limit point $P_L$. Let $P\neq P_L$ be a point of $\mathcal{C}$ and let $H_P(t)$ be its fibre function (under $\mathcal{P}(t)$). Then, $\mult_P(\mathcal{C})=\deg(H_P(t)).$
\end{corollary}

\para

\begin{example} Let $\cal C$ be the rational plane curve defined by the projective parametrization
$$\mathcal{P}(t)=(-t^3-5t^2-7t-3:t^4+7t^3+17t^2+17t+6:t^4+1)\in\mathbb{P}^2(\mathbb{C}(t)).$$
Let us compute $H_P(t)$, where   $P=(0:0:1)$.
Since $P$ is an affine point, we can obtain  $H_P(t)$ from $\phi_1$ and $\phi_2$ (see Remark \ref{R-afin-noafin}). Since $\phi_1(t)=-2 t^3-10 t^2-14 t-6$ and $\phi_2(t)=2 t^4+14
t^3+34 t^2+34 t+12$, we get that
$$H_P(t)=\gcd(\phi_1,\phi_2)=2(t+3)(t+1)^2.$$

\begin{figure}[h]
$$
\begin{array}{cc}
\psfig{figure=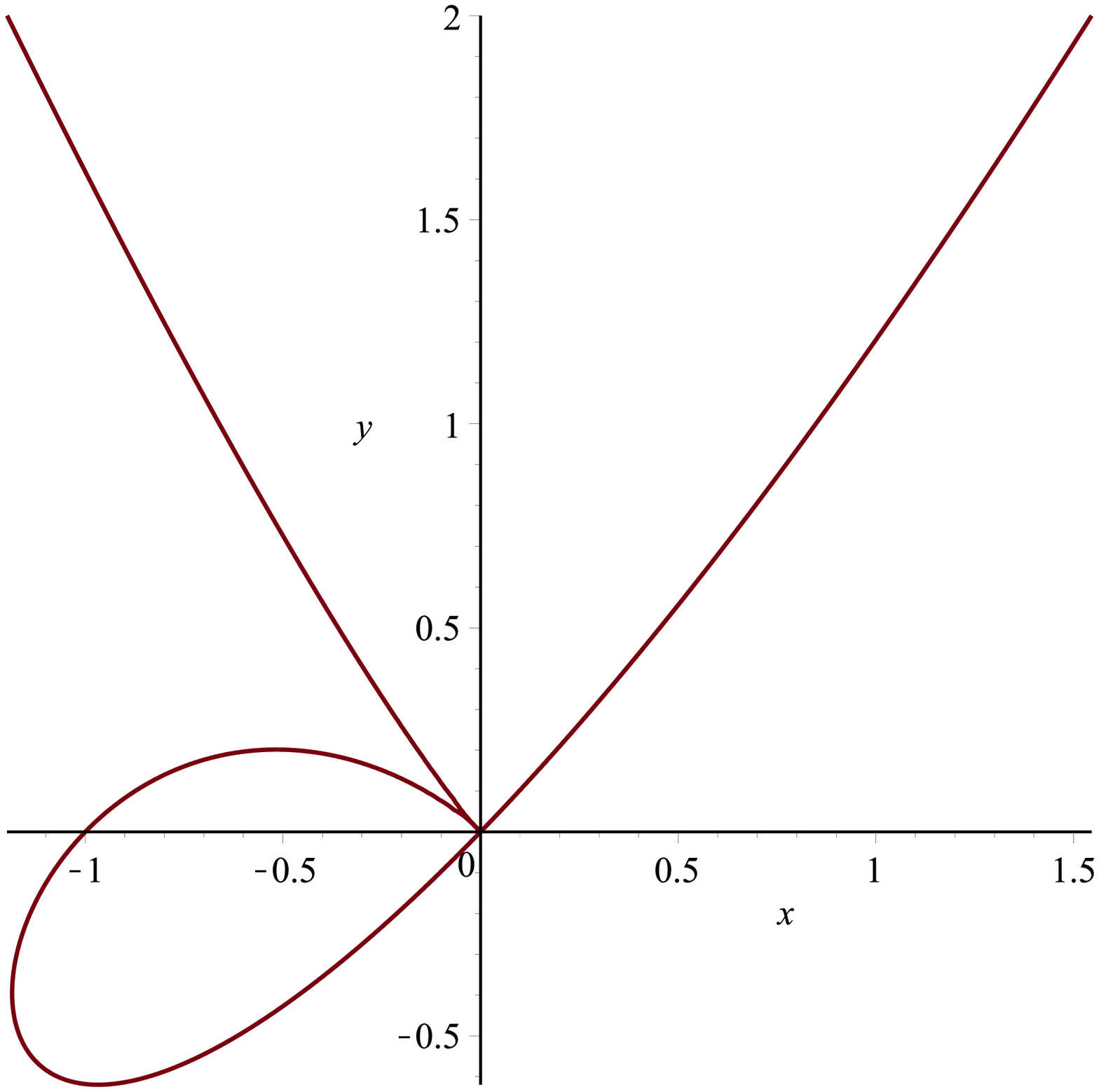,width=4.5cm,height=4.5cm,angle=0} &
\psfig{figure=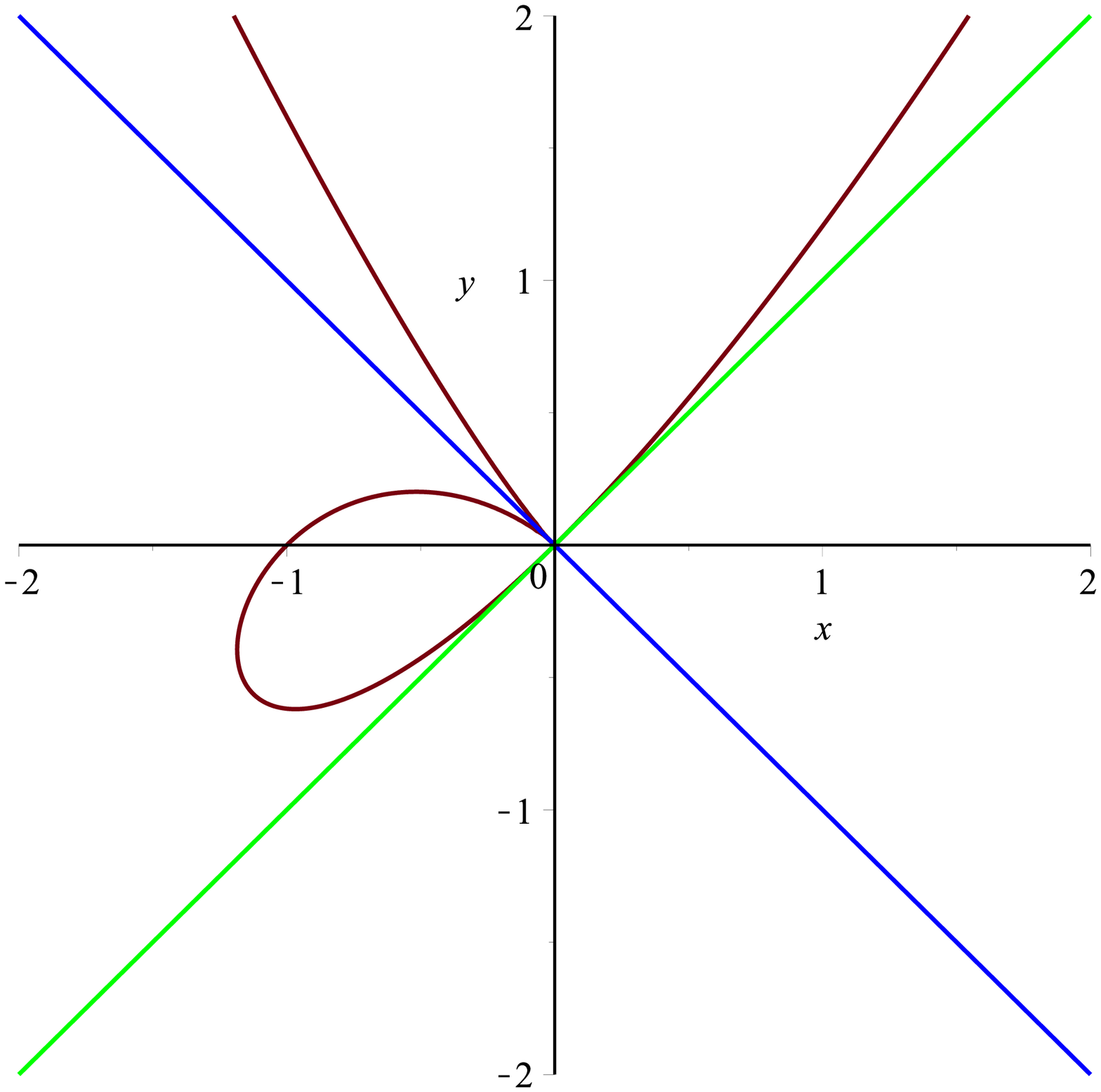,width=4.5cm,height=4.5cm,angle=0}
\end{array}
$$ \caption{Triple point with two tangents}\label{F-ptotriplenoord}
\end{figure}

Therefore, $\cP(-1)=\cP(-3)=P$, and applying Theorem \ref{T-tangentes}, we deduce that $\cal C$ has at $P$ two different tangents, one of multiplicity $1$ and the other one of multiplicity $2$. The parametrizations defining these tangents are given as $$\tau_1(t)=\mathcal{P}(-3)+\mathcal{P}'(-3)t\quad\text{ and}\quad \tau_2(t)=\mathcal{P}(-1)+\frac{\mathcal{P}''(-1)}{2}t^2,$$
respectively (see \cite{JSC-Perez}). Note that these tangents are the lines  $y=x$ and $y=-x$ (see Figure \ref{F-ptotriplenoord}). Finally, we conclude that  $P$ is a non--ordinary point of multiplicity $3$ (see Corollary \ref{C-multiplicidad}).

\end{example}

%In Theorem \ref{T-tangentes} and Corollary \ref{C-multiplicidad} we assume that $P\neq P_L$; otherwise, an additional tangent may arise which is not associated to any value of the fibre. This circumstance will be analyzed in %the next section.

\section{Resultants and singularities}\label{S-resultantsingularities}

In Section \ref{S-notacion}, we show that, given a rational proper parametrization,
$\mathcal{P}(t)$, the multiplicity of a given point, $P\neq P_L$,  is the cardinality of the fibre of  $\mathcal{P}(t)$  at $P$ (see Corollary
\ref{C-multiplicidad}). That is,  the multiplicity of $P=\mathcal{P}(s_0),\,s_0\in\mathbb{K}$  is given by the cardinality of the set
$$\mathcal{F}_{\mathcal{P}}(\mathcal{P}(s_0))=\{t\in\mathbb{K}:\mathcal{P}(t)=\mathcal{P}(s_0)\} $$
(note that we are assuming  that $P\neq P_L$). Observe that $s_0\in \mathcal{F}_{\mathcal{P}}(\mathcal{P}(s_0))$ and hence, the cardinality of $\mathcal{F}_{\mathcal{P}}(\mathcal{P}(s_0))$ is greater than or equal to $1$. Thus, $\mathcal{P}(s_0)$ is a singular point if and only if  the cardinality of $\mathcal{F}_{\mathcal{P}}(\mathcal{P}(s_0))$ is greater than $1$.

\para

Taking into account the above statement, in this section, we show how the different factors of a univariate resultant computed from the polynomials $G_i(s,t),\,i=1,2,3,$ are exactly the fibre functions of the singularities  of $\cal C$. Thus, in particular, the singularities  of $\cal C$ and its corresponding multiplicities are determined.   The idea for the construction of the resultant is that a point $\mathcal{P}(s_0)\in {\cal C},\,s_0\in\mathbb{K},$ is a singularity if and only if $\deg(H_{\mathcal{P}(s_0)}(t))>1$ (i.e. the fibre equations of $\mathcal{P}(s_0)$ have more than one
common solution).

\para

For this purpose, we first assume that $\mathcal{P}(s_0)$ is an affine point. Thus, Remark \ref{R-afin-noafin} implies that the fibre equations are given by
$$\left\{\begin{array}{l}p_1(t)p(s_0)-p_1(s_0)p(t)=0\\p_2(t)p(s_0)-p_2(s_0)p(t)=0.\end{array}\right.$$
Note that this is equivalent to $G_1(s_0,t)=G_2(s_0,t)=0$, where  $G_1(s,t)$ and $G_2(s,t)$ are the polynomials introduced in
(\ref{Eq-fibra-generica}).

\para

Then, $\mathcal{P}(s_0)$ is a singular point if and only if $G_1(s_0,t)$ and $G_2(s_0,t)$ have more than one common root or, equivalently, if and only if the polynomials  $G_1(s_0,t)/(t-s_0)$ and $G_2(s_0,t)/(t-s_0)$ have a common root (we note that $s_0$ is already a root of $G_1(s_0,t)$ and $G_2(s_0,t)$). This implies  that
$$\Res_t\left(\frac{G_1(s_0,t)}{t-s_0},\frac{G_2(s_0,t)}{t-s_0}\right)=0.$$
Hence, given the polynomial
$$R(s)=\Res_t\left(\frac{G_1(s,t)}{t-s},\frac{G_2(s,t)}{t-s}\right),$$
if the point $\mathcal{P}(s_0)$ is singular, then $R(s_0)=0$. In fact, in \cite{Abhy}, it is proved that this resultant provides the product of the fibre functions of the singularities of the curve, in the case that ${\cal P}(t)$ is a polynomial parametrization. A  generalization for the case of a given rational parametrization  (not necessarily polynomial) is presented in \cite{Buse2012}.

%\para
%Recall that we are assuming $\mathcal{P}(t)$ to be
%a proper parametrization. Otherwise, the singular points are those where the cardinality of the fibre is greater than the tracing index. In this case (see \cite{JSC-Perez}) $R(s)$ must be replaced by
%$$T(s):=\Res_t\left(\frac{G_1(s,t)}{G(s,t)},\frac{G_2(s,t)}{G(s,t)}\right),$$
%where $G(s,t)=\gcd(G_1(s,t),G_2(s,t))$.

\para

Thus,  $R(s)$ can be used to compute the singularities of the curve, but some problems could appear. First, the values $s\in\mathbb{K}$ that provide singular points are roots of the polynomial $R$ but the reciprocal is not true; i.e. a root of $R$ may not provide a singular point. In addition, we are assuming that the singularity is an affine point, but also singularities at infinity have to be detected.

\para

 The {\it  T--function}, that we introduce below, improves the properties of $R(s)$ and characterizes the singular points of $\cal C$ (affine and at infinity). In order to introduce it, we need to consider  $$\delta_i:=\degree_t(G_i),\quad
\lambda_{ij}:=\min\{\delta_i,\delta_j\},\quad
G_i^*(s,t):=\displaystyle\frac{G_i(s,t)}{t-s}\in {\Bbb K}[s,t]$$
and
$$R_{ij}(s):=\Res_t(G_i^*,G_j^*)\in {\Bbb K}[s]\,\, \, \, \mbox{for}\,\,\, i,j=1,2,3,\,\,i<j.$$

\begin{definition}\label{D-T-funct}
We define the T--function of the parametrization  $\mathcal{P}(t)$
as
$$T(s)=R_{12}(s)/p(s)^{\lambda_{12}-1}.$$
\end{definition}

\para

In the following we show that this function provides essential information concerning the singularities of the given curve $\cal C$ (see Theorem \ref{T-tfunct}). To start with, the following proposition claims that the T--function can be defined similarly from  $R_{13}(s)$ or $R_{23}(s)$. In addition, in  Corollary \ref{C-T-polin}, we prove that $T(s)$ is a polynomial.

%In fact, as we stated above, the factorization of the T--function provides the fibre functions of each singularity as well as its corresponding multiplicity. We remark that from the fibre function of a point $P$, one may determine the multiplicity of $P$,  its fibre ${\mathcal F}_{\mathcal P}(P)$, and the tangent lines of $\cal C$ at $P$ (see Section \ref{S-notacion}).

\para

\begin{proposition}\label{P-res13-res23}
It holds that
$$T(s)=\frac{R_{12}(s)}{p(s)^{\lambda_{12}-1}}=\frac{R_{13}(s)}{p_1(s)^{\lambda_{13}-1}}=\frac{R_{23}(s)}{p_2(s)^{\lambda_{23}-1}}.$$
\end{proposition}

\noindent\textbf{Proof:} We distinguish two steps to prove the proposition:

\begin{center}
{\it Step 1}
\end{center}

\noindent First, we show that
$$\frac{R_{12}(s)}{p(s)^{\lambda_{12}-1}}=\frac{R_{13}(s)}{p_1(s)^{\lambda_{13}-1}}.$$
For this purpose, we see the polynomial $G_1(s,t)\in {\Bbb K}[s,t]$ as a polynomial in the variable $t$ that is, $G_1(s,t)\in ({\Bbb K}[s])[t]$. Since $\deg_t(G_1)=\deg_s(G_1)=\delta_1$  (see Remark \ref{R-grado-Gs}), then $G_1$ has $\delta_1$ roots (in the variable $t$), and one of them is $t=s$. Thus, we may write
$$G_1(s,t)=\lc_t(G_1)(t-s)(t-\alpha_1(s))\cdots
(t-\alpha_{\delta_1-1}(s))$$ and
\begin{equation}\label{Eq-factorG1}G^*_1(s,t)=\lc_t(G_1)(t-\alpha_1(s))\cdots
(t-\alpha_{\delta_1-1}(s)),\end{equation}
where $\lc_t(\cdot)$ denotes the leader coefficient with respect to the variable $t$ of a polynomial $(\cdot)$. Now, taking into account the properties of the resultants (see e.g. \cite{Cox1998}, \cite{SWP}, \cite{Vander}), we get that
\begin{equation}\label{Eq-res12}R_{12}(s):=\Res_t(G_1^*,G_2^*)=\lc_t(G^*_1)^{\delta_2-1}\prod_{i=1}^{\delta_1-1}G^*_2(s,\alpha_i(s)).\end{equation}
Note that $G_2(s,t)=p_2(s)p(t)-p(s)p_2(t)$, and thus
$G_2(s,\alpha_i(s))=p_2(s)p(\alpha_i(s))-p(s)p_2(\alpha_i(s)).$ Furthermore, since
$G_1(s,\alpha_i(s))=p_1(s)p(\alpha_i(s))-p(s)p_1(\alpha_i(s))=0$,
we get that
$$p(\alpha_i(s))=\frac{p_1(\alpha_i(s))}{p_1(s)}p(s).$$
Therefore, $$
G_2(s,\alpha_i(s))=p_2(s)\frac{p_1(\alpha_i(s))}{p_1(s)}p(s)-p(s)p_2(\alpha_i(s))=$$$$(p_1(\alpha_i(s))p_2(s)-p_2(\alpha_i(s))p_1(s))\frac{p(s)}{p_1(s)}=G_3(s,\alpha_i(s))\frac{p(s)}{p_1(s)}$$
which implies that
$$G^*_2(s,\alpha_i(s))=G^*_3(s,\alpha_i(s))\frac{p(s)}{p_1(s)}.$$
Now we substitute in (\ref{Eq-res12}) and we get
$$R_{12}(s)=\left(\lc_t(G^*_1)^{\delta_2-1}\prod_{i=1}^{\delta_1-1}G^*_3(s,\alpha_i(s))\right)\left(\frac{p(s)}{p_1(s)}\right)^{\delta_1-1},$$
which can be expressed as
%$$R_{12}(s)=\left(lc_t(G^*_1)^{\delta_3-1}\prod_{i=1}^{\delta_1-1}G^*_3(s,\alpha_i(s))\right)lc_t(G^*_1)^{\delta_2-\delta_3}\left(\frac{p(s)}{p_1(s)}\right)^{\delta_1-1}$$
%y, por tanto,
$$R_{12}(s)=R_{13}(s)\lc_t(G^*_1)^{\delta_2-\delta_3}\left(\frac{p(s)}{p_1(s)}\right)^{\delta_1-1}.$$
Hence, we only have to prove that
\begin{equation}\label{Eq-aux1}\lc_t(G^*_1)^{\delta_2-\delta_3}\left(\frac{p(s)}{p_1(s)}\right)^{\delta_1-1}=\frac{p(s)^{\lambda_{12}-1}}{p_1(s)^{\lambda_{13}-1}},\end{equation}
and, for this purpose, we consider different cases depending on   $d_1$, $d_2$ and $d_3$ (that is, on the degrees of $p_1$,
$p_2$ and $p$). We  remind that  $\delta_1=\max\{d_1,d_3\}$,
$\delta_2=\max\{d_2,d_3\}$ and $\delta_3=\max\{d_1,d_2\}$ (see
Remark \ref{R-grado-Gs}).

\begin{itemize}
\item {Case 1:} Let $d_1<d_3$. Then,
$\delta_1=d_3\leq \delta_2$, and $\lambda_{12}=\delta_1$.
In addition, it holds that $\lc_t(G^*_1)=\lc_t(G_1)=p_1(s)$ since
$G_1(s,t)=p_1(s)p(t)-p(s)p_1(t)$ and $d_1<d_3$. Thus,
$$\lc_t(G^*_1)^{\delta_2-\delta_3}\left(\frac{p(s)}{p_1(s)}\right)^{\delta_1-1}=\frac{p(s)^{\delta_1-1}}{p_1(s)^{\delta_1-\delta_2+\delta_3-1}}
=\frac{p(s)^{\lambda_{12}-1}}{p_1(s)^{\delta_1-\delta_2+\delta_3-1}}.$$
In order to check that (\ref{Eq-aux1}) holds, we only have to prove that
$\delta_1-\delta_2+\delta_3=\lambda_{13}$. Let us see that this equality holds in the following situations:
\begin{itemize}
\item[a)] $\delta_1<\delta_3$: then,  $d_2>d_1,d_3$ which implies that $\delta_2=\delta_3=d_2$ and $\lambda_{13}=\delta_1$.
\item[b)] $\delta_1>\delta_3$: then,  $d_3>d_1,d_2$ which implies that $
\delta_1=\delta_2=d_3$ and $\lambda_{13}=\delta_3$.
\item[c)] $\delta_1=\delta_3$: then,  $d_3=\max\{d_1,d_2\}$ which implies that $
d_1<d_2=d_3$ and $\delta_1=\delta_2=\delta_3$.
\end{itemize}
\item {Case 2:} Let $d_1>d_3$. Then,  $\delta_1=d_1\leq \delta_3$, which implies that $\lambda_{13}=\delta_1$. In addition,
$\lc_t(G^*_1)=p(s)$, and then
$$\lc_t(G^*_1)^{\delta_2-\delta_3}\left(\frac{p(s)}{p_1(s)}\right)^{\delta_1-1}=\frac{p(s)^{\delta_1+\delta_2-\delta_3-1}}{p_1(s)^{\delta_1-1}}
=\frac{p(s)^{\delta_1+\delta_2-\delta_3-1}}{p_1(s)^{\lambda_{13}-1}}.$$
Thus, we only have to prove that
$\delta_1+\delta_2-\delta_3=\lambda_{12}$. For this purpose, we reason similarly as in Case 1 by considering the following situations:
$\delta_1<\delta_2$, $\delta_1>\delta_2$ and $\delta_1=\delta_2$.
\item {Case 3:} Let $d_1=d_3<d_2$.  Then,
$\delta_2=\delta_3=d_2$, and thus $\lc_t(G^*_1)^{\delta_2-\delta_3}=1$.
In addition, $\delta_1\leq\delta_2=\delta_3$, which implies that
$\lambda_{12}=\lambda_{13}=\delta_1$.
\item { Case 4:} Let $d_2<d_1=d_3$. In this case, we have that
$\delta_1=\delta_2=\delta_3$ and then, (\ref{Eq-aux1}) trivially holds.
\item {Case 5:} Let $d_1=d_2=d_3$. This case is similar to Case 4.
\end{itemize}

\begin{center}
{\it Step 2}
\end{center}

\noindent Let us prove that
$$\frac{R_{12}(s)}{p(s)^{\lambda_{12}-1}}=\frac{R_{23}(s)}{p_2(s)^{\lambda_{23}-1}}.$$
For this purpose, we observe that, up to constants in ${\Bbb K}\setminus\{0\}$, it holds that
$R_{12}(s)=R_{21}(s)$. Thus, we may write
$$R_{12}(s)=R_{21}(s)=\lc_t(G^*_2)^{\delta_1-1}\prod_{i=1}^{\delta_2-1}G^*_1(s,\beta_i(s))$$
where
$$G^*_2(s,t)=\lc_t(G_2)(t-\beta_1(s))\cdots
(t-\beta_{\delta_2-1}(s)).$$ Now, we observe that these equalities are equivalent to   (\ref{Eq-res12}) and (\ref{Eq-factorG1}), respectively. Thus, reasoning similarly as above, we obtain that
$$R_{12}(s)=R_{23}(s)\lc_t(G^*_2)^{\delta_1-\delta_3}\left(\frac{p(s)}{p_2(s)}\right)^{\delta_2-1}$$
and that
$$\lc_t(G^*_2)^{\delta_1-\delta_3}\left(\frac{p(s)}{p_2(s)}\right)^{\delta_2-1}=\frac{p(s)^{\lambda_{12}-1}}{p_2(s)^{\lambda_{23}-1}}.$$
\hfill $\Box$

\para

\begin{corollary}\label{C-T-polin} It holds that
$T(s)\in {\Bbb K}[s]$.
\end{corollary}

\noindent\textbf{Proof:} Let us assume that $T(s)$ is not a polynomial. Then, we simplify the rational function and we write
$$\frac{R_{12}(s)}{p(s)^{\lambda_{12}-1}}=\frac{M_{12}(s)}{\bar{p}(s)},$$ where  $M_{12}(s)\in {\Bbb K}[s]$, $\overline{p}(s)\in {\Bbb K}[s]\setminus{\Bbb K}$ and  $\gcd(M_{12},\bar{p})=1$. Note that  $\overline{p}$ divides $p^{\lambda_{12}-1}$.
%which implies that lo que implica que los distintos factores de $\overline{p}$ deben dividir a $p$.

\para

\noindent Similarly, from Proposition \ref{P-res13-res23}, we have that
$$\left\{\begin{array}{cc}
\displaystyle\frac{R_{13}(s)}{p_1(s)^{\lambda_{13}-1}}=\frac{M_{13}(s)}{\bar{p}_1(s)} & \text{where } \overline{p}_1 \text{ divides } p_1^{\lambda_{13}-1}, \text{ and }\gcd(M_{13},\bar{p}_1)=1 \\
\displaystyle\frac{R_{23}(s)}{p_2(s)^{\lambda_{23}-1}}=\frac{M_{23}(s)}{\bar{p}_2(s)} & \text{where } \overline{p}_2 \text{ divides } p_2^{\lambda_{23}-1}, \text{ and }\gcd(M_{23},\bar{p}_2)=1.
\end{array}\right.
$$
Furthermore, we have that (see Proposition \ref{P-res13-res23})
$$\frac{M_{12}(s)}{\bar{p}(s)}=\frac{M_{13}(s)}{\bar{p}_1(s)}=\frac{M_{23}(s)}{\bar{p}_2(s)}$$ which implies that
 $$M_{12}(s)\bar{p}_1(s)=M_{13}\bar{p}(s) \mbox{ and } M_{23}(s)\bar{p}_1(s)=M_{13}\bar{p}_2(s).$$ Taking into account that $\gcd(M_{12},\bar{p})=\gcd(M_{13},\bar{p}_1)=\gcd(M_{23},\bar{p}_2)=1$, and the above equalities, we get that $\overline{p}_1=\overline{p}_2=\overline{p}$. Then, we deduce that   $\overline{p}$ divides $p$, $p_1$ and $p_2$, which is impossible since  $\gcd(p_1,p_2,p)=1$.
\hfill $\Box$

\para

In the following theorem, we show how the ordinary singularities of $\cal C$ can be determined from the T--function. In fact, $T(s)$ describes totally the singularities of the curve, since its factorization provides the fibre functions of each singularity as well as its corresponding multiplicity.  From the fibre function, $H_P(t)$, of a point $P$, one obtains the multiplicity of $P$, its fibre,  and the tangent lines at $P$ (see Section \ref{S-notacion}).

\para

An alternative approach for computing this factorization, based on the construction of $\mu$--basis, can be found in \cite{Buse2012}.

\para

In Theorem \ref{T-tfunct}, we assume that $\cal C$  has only ordinary singularities. Otherwise, for applying this theorem, we should apply quadratic transformations (blow-ups) for birationally transforming $\cal C$ into a curve with only ordinary singularities (see Chapter 2 in \cite{SWP}). For such a curve the following theorem holds.

\para

\begin{theorem}\label{T-tfunct} {\bf(Main theorem)} Let $\mathcal{C}$ be a rational algebraic curve defined by a parametrization $\mathcal{P}(t)$, with limit point $P_L$. Let $P_1,\ldots,P_n$ be the singular points of $\cal C$, with multiplicities $m_1,\ldots,m_n$ respectively. Let us assume that they are ordinary singularities and that $P_i\neq P_L$ for $i=1,\ldots,n$. Then, it holds that
$$T(s)=\prod_{i=1}^nH_{P_i}(s)^{m_i-1}.$$
\end{theorem}

\para

This theorem will be proved in Section  \ref{S-proof} and a generalization for the case of space curves of any dimension will be presented in Section \ref{S-spacialcurves}. Moreover, in \cite{MyB-2017(b)}, we prove that the theorem holds also if $P_L$ is a singularity. Finally, an analogous result which admits the existence of non--ordinary singularities in the curve will be developed in a future work.

\para

\begin{corollary}\label{C-grado-T}
Let $\cal C$  be a rational plane curve such that all its singularities are ordinary. Let $\mathcal{P}(t)$ be a parametrization of $\cal C$ such that $P_L$ is regular. It holds that
$\deg(T)=(d-1)(d-2).$
\end{corollary}

\noindent\textbf{Proof:} From Theorem \ref{T-tfunct} and Corollary  \ref{C-multiplicidad}, we have that  $\deg(T)=\sum_{i=1}^nm_i(m_i-1),$ where $P_1,\ldots,P_n$
are the singular points, and  $m_1,\ldots,m_n$ its corresponding multiplicities. Since  $\mathcal{C}$ is a rational curve, its genus is $0$ and thus $\sum_{i=1}^nm_i(m_i-1)=(d-1)(d-2)$ (see Chapter 3 in \cite{SWP}). \hfill $\Box$

\para

%In the following example, Theorem \ref{T-tfunct} proves to be really useful for getting information about the singularities of a rational curve.

\para

\begin{example}\label{E-Tfunct-planas}
Let $\mathcal{C}$ be the rational plane curve defined by the projective parametrization
$$\mathcal{P}(t)=(t^5-5t^4+5t^3+5t^2-6t:t^5+5t^4+5t^3-5t^2-6t:t^4-13t^2+36)\in\mathbb{P}^2(\mathbb{C}(t)).$$
We compute the T--function and we get that, up to constants in $\mathbb{C}$,
$$T(s)=(t-2)(t-3)(t+3)(t+2)t^2(t^2+6)(t-1)^2(t+1)^2.$$
Thus, from Theorem \ref{T-tfunct}, we deduce that the fibre function of each singularity of  $\mathcal{C}$  appears in the polynomial $T$. Let us analyze the different factors of $T$:

\begin{itemize}
\item The factors with power $2$ correspond to triple points. Indeed: these factors are $t$, $(t-1)$ and $(t+1)$, and we have that $P_1=\mathcal{P}(0)=\mathcal{P}(1)=\mathcal{P}(-1)=(0:0:1)$. Then, $P_1$  is a triple point whose fibre function is   $$H_{P_1}(t)=(t-1)(t+1)t.$$
\item  The factors with power $1$ correspond to different double points. In order to determine the associated factors, we should compute the corresponding fibre functions. For  the factors $(t-2)$ and $(t-3)$, we have that $P_2=\mathcal{P}(2)=\mathcal{P}(3)=(0:1:0)$ and thus $$H_{P_2}(t)=(t-2)(t-3).$$
   This implies that  $P_2$ is a double point at  infinity. %and it is reached by  $\mathcal{P}(2)$ and $\mathcal{P}(3)$.
\item Similarly, if we consider the factors  $(t+2)$ and  $(t+3)$, we get that  $P_3=\mathcal{P}(-2)=\mathcal{P}(-3)=(1:0:0)$, and its fibre function is $$H_{P_3}(t)=(t+2)(t+3),$$
   which implies that $P_3$ is a double point at infinity. % and it is reached by  $\mathcal{P}(-2)$ and $\mathcal{P}(-3)$.
\item Finally, the factor $(t^2+6)$  provides the point $P_4=\mathcal{P}(-I\sqrt{6})=\mathcal{P}(+I\sqrt{6})=(-7/5:7/5:1)$. Hence, $P_4$ is an affine double point and its fibre function is $$H_{P_4}(t)=(t^2+6).$$
\end{itemize}

\begin{figure}[h]
$$
\begin{array}{cc}
\psfig{figure=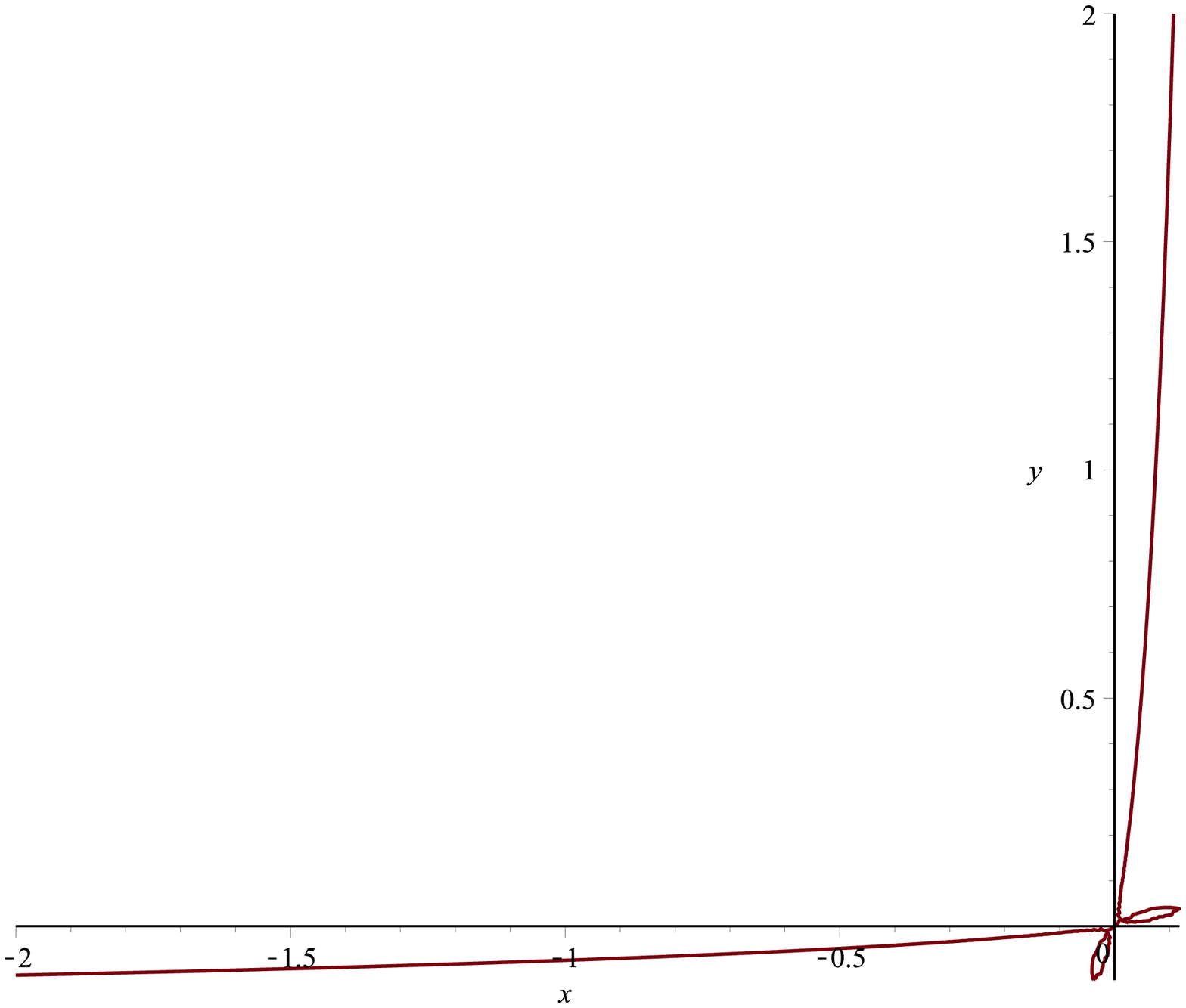,width=6cm,height=5cm,angle=0} &
\psfig{figure=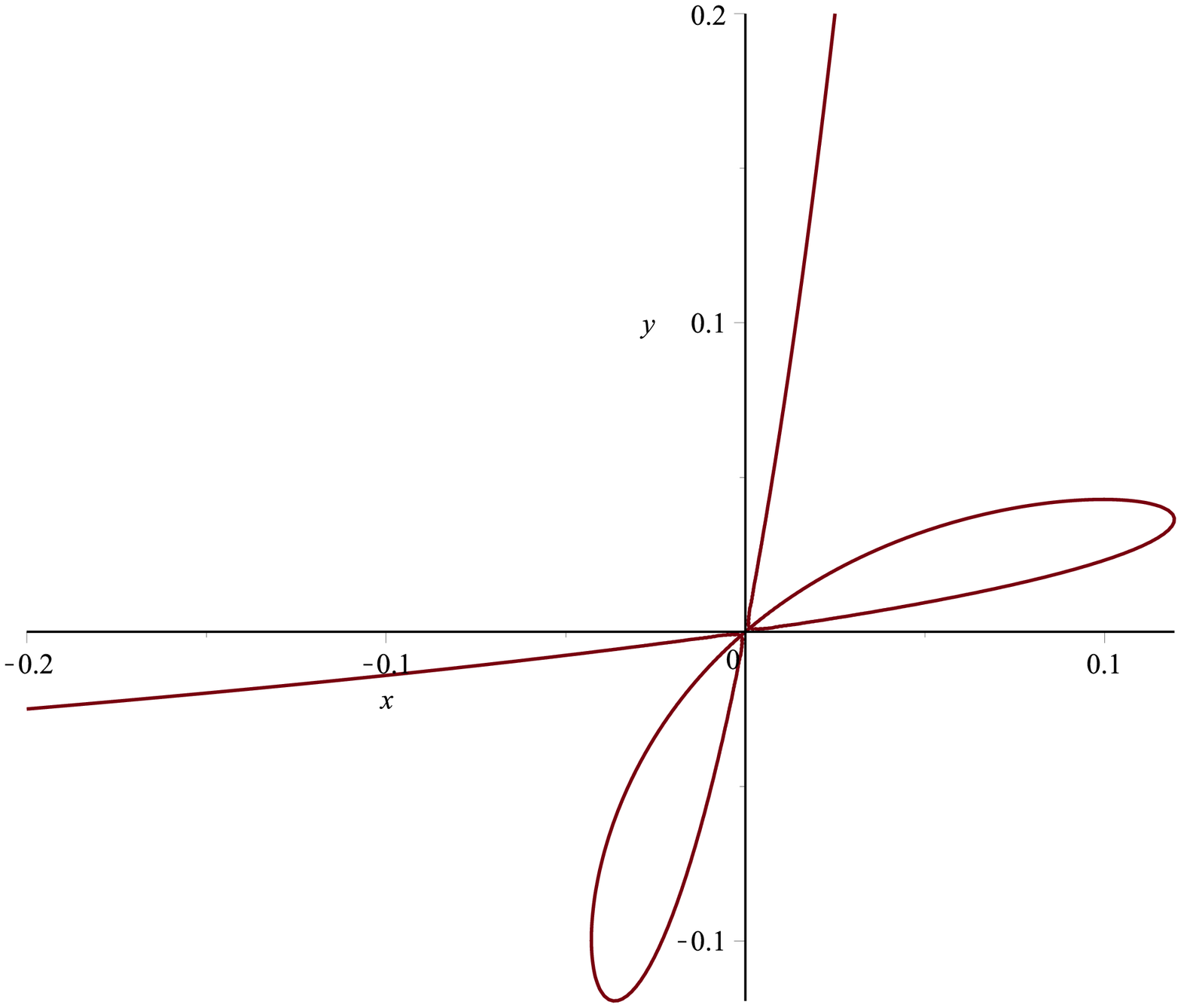,width=6cm,height=5cm,angle=0}
\end{array}
$$ \caption{Curve $\cal C$ (left) and a neighborhood of the triple point $P_1$ (right)}\label{F-T1}
\end{figure}

\noindent Note that,
$$T(s)=H_{P_1}(s)^{m_1-1} H_{P_2}(s)^{m_2-1} H_{P_3}(s)^{m_3-1} H_{P_4}(s)^{m_4-1}.$$
Furthermore, we observe that Corollary \ref{C-grado-T} is verified. Indeed, we have that $d=5$ and $\degree_t(T)=12=(d-1)(d-2)$.\\

In Figure \ref{F-T1}, we plot the curve  $\mathcal{C}$, and a neighborhood of the triple point $P_1$. Observe that $P_4$ is an isolated point.

\end{example}

\para

In Example \ref{E-Tfunct-planas}, we have been able to determine the singularities of $\cal C$ and its corresponding multiplicities, by computing the factors of the polynomial $T(s)$. However, in  general, one needs to
introduce algebraic numbers during the computations. In the
following, we present a method that allows us to determine the factors of the polynomial $T(s)$ (and thus, the  singularities of a curve) without directly
using algebraic numbers. For this purpose, we introduce the notion of  {\it family of conjugate parametric
points} (see \cite{JSC-Perez}), which generalizes the concept of family of conjugate points (see e.g.  Chapter 2 in \cite{SWP}). The  idea is to collect points whose coordinates depend
algebraically on all the conjugate roots of the same irreducible polynomial $m(t)$. The computations on such a family of points can be carried out by using the polynomial $m(t)$.

\para

\begin{definition}\label{D-conjugatepoints}
Let
$$\mathcal{G}=\{(p_1(\alpha):p_2(\alpha):p(\alpha))|m(\alpha)=0\}\subset \mathbb{P}^2(\mathbb{K}).$$
The set $\mathcal{G}$ is called a family of conjugate parametric
points over $\mathbb{K}$ if the following conditions are satisfied:
\begin{enumerate}
\item $p_1,p_2,p,m\in\mathbb{K}[t]$ and $\gcd(p_1,p_2,p)=1$.
\item $m$ is irreducible.
\item $\deg(p_1),\deg(p_2),\deg(p)<\deg(m)$.
\end{enumerate}
We denote such a family by $\mathcal{G}=\{\mathcal{P}(s)\}_{m(s)}=\{(p_1(s):p_2(s):p(s))\}_{m(s)}.$
\end{definition}

\para

Condition $2$ in Definition \ref{D-conjugatepoints} can be  stated considering that $m$ is only square-free (see  Definition 12 in  \cite{JSC-Perez}). However, in order to prove Theorem \ref{T-characterconjugate}, one needs  $m$ to be an irreducible polynomial (see Theorem 16 in \cite{JSC-Perez}). Hence, using the above definition, in \cite{JSC-Perez} it is proved the following theorem.

\para

\begin{theorem}\label{T-characterconjugate}
The singularities of the  curve $\cal C$ can be decomposed as a
finite union of families of conjugate parametric points over
$\mathbb{K}$ such that all points in the same family have the same
multiplicity and character.
\end{theorem}

\para

If some singularities of the given curve are in a family
$\mathcal{G}=\{\mathcal{P}(s)\}_{m(s)}=\{(p_1(s):p_2(s):p(s))\}_{m(s)}$, then the polynomial $m(s)$ will be an irreducible factor of the T--function. In this case, Theorem \ref{T-ptos-conjug} allows us to determine the singularities provided by $\mathcal{G}$ and their corresponding multiplicities.

\para

\begin{theorem}\label{T-ptos-conjug}
Let $m(s)$ be an irreducible polynomial such that   $m(s)^{k-1}, k\in {\Bbb N},\,k\geq 1,$ divides $T(s)$. Then, the roots of $m(s)$ determine the fibre of some singularities of multiplicity $k$ that are defined by a  family of conjugate parametric points. The number of singularities in such a family is  $n=\deg(m(s))/k$.
\end{theorem}

\noindent\textbf{Proof:} From Theorem \ref{T-tfunct}, we get that the points in $\mathcal{G}=\{\mathcal{P}(s)\}_{m(s)}$ are singularities of multiplicity  $k$. In addition, if $\mathcal{G}=\{P_1,\ldots,P_n\}$, we have that
$$m(s)=\prod_{i=1}^nH_{P_i}(s),$$
where $H_{P_1},\ldots,H_{P_n}$ are the fibre functions of the points $P_1, \ldots, P_n,$ respectively. From Corollary  \ref{C-multiplicidad}, we get that $\deg(H_{P_i})=\mult_{P_i}(\mathcal{C}),\,i=1,\ldots,n,$ and since, in this case, $\mult_{P_i}(\mathcal{C})=k$, we conclude   that $\deg(m(s))=nk.$ \hfill $\Box$

\para

\begin{example}
Let $\cal C$ be the rational curve defined by $\mathcal{P}(t)=(p_1(t):p_2(t):p(t))\in\mathbb{P}^2(\mathbb{C}(t))$, where
$$\begin{array}{l}p_1(t)=t^5-13t^4+63t^3-143t^2+152t-60\\
p_2(t)=t^5-21t^4+157t^3-507t^2+706t-336\\
p(t)=t^5+7t^4-2t^3+t-1.\end{array}$$ The T--function is given by\\

\noindent
$T(s)=28161216(968585964-2319881360s+2070988203s^2-904722886s^3+208513387s^4-24407436s^5+1145528s^6)(s-1)^2(s-2)^2(s-3)^2.$\\

\noindent From the polynomial $T(s)$, we deduce that the singularities of $\cal C$ are:

\begin{itemize}
\item The triple point
$P_1=\mathcal{P}(1)=\mathcal{P}(2)=\mathcal{P}(3)=(0:0:1)$.
The fibre function of $P_1$ is  $H_{P_1}(s)=(s-1)(s-2)(s-3)$, and these factors appear with power $2$ in the polynomial $T(s)$.
\item Three double points associated to the irreducible factor\\
\noindent
$m(s)=968585964-2319881360s+2070988203s^2-904722886s^3+208513387s^4-24407436s^5+1145528s^6.$\\
Since $m(s)$ appears with power  $1$ in the polynomial $T(s)$, we conclude that it is associated to double points. In addition, using Theorem \ref{T-ptos-conjug}, we get that this factor provides three different points  (each point has a fibre function of degree $2$ and the three fibre functions have power $1$; multiplying these fibre functions we obtain the polynomial $m$).
%Ahora, dividiendo su grado (6) por la multiplicidad (2), deducimos que contiene 3 puntos
\end{itemize}
In Figure \ref{F-1x3+3x2}, we plot the given curve $\cal C$, and we can see the singularities obtained   (note that each singularity is real and affine).

\begin{figure}[h]
$$
\psfig{figure=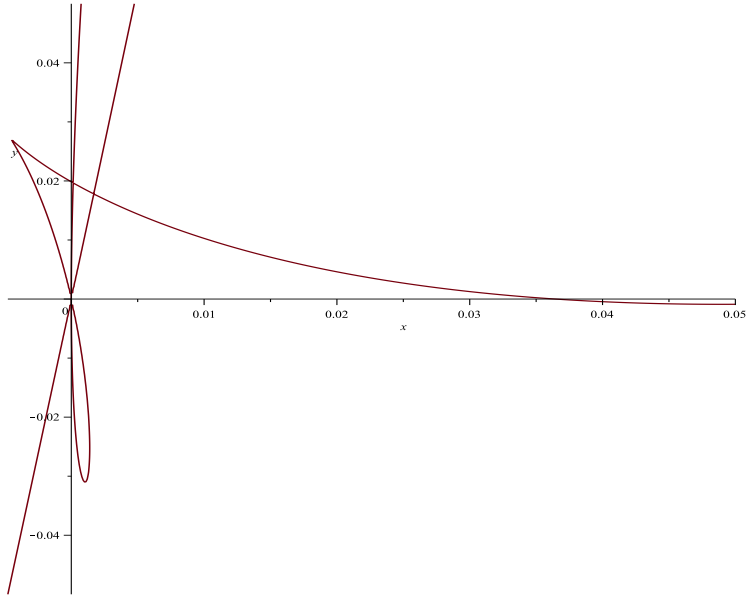,width=10cm,height=8cm}
$$ \caption{The curve $\cal C$ that has a triple point and three double points}\label{F-1x3+3x2}
\end{figure}
\end{example}

\section{The general case for  rational space curves}\label{S-spacialcurves}

In this section, we show that Theorem \ref{T-tfunct} can be applied for the case that the given curve $\cal C$  is a rational space curve. In this case, we construct an equivalent polynomial to the T--function  introduced for plane curves (see Definition \ref{D-T-funct}),  and we prove that this polynomial, which will be denoted as $T_E(s)$,  describes totally the singularities of $\cal C$,   since each factor of  $T_E(s)$  is a power of the fibre function of one singularity of the given curve. This power is, in fact, the multiplicity of the singularity minus 1. We recall that from the fibre function of a point $P$, one may determine the multiplicity of $P$  as well as its fibre  ${\mathcal F}_{\mathcal P}(P)$  and the tangent lines of $\cal C$ at $P$ (see Section \ref{S-notacion}). The method presented generalizes the results obtained in \cite{Rubio}, since a complete classification of the singularities of a given space curve, via the factorization of a univariate  resultant, is obtained.

\para

\noindent
In the following, we consider
$$\mathcal{P}(t)=(p_1(t):\cdots: p_n(t):p(t))\in\mathbb{P}^n(\mathbb{K}(t)),\quad \gcd(p_1,\ldots,p_n,p)=1,$$
a proper parametrization of a given rational space curve ${\cal C}$. In addition, we define the associated rational parametrization over ${\Bbb K}(Z)$, where $Z=(Z_1,\ldots,Z_{n-2})$ and $Z_1,\ldots,Z_{n-2}$ are new variables, given by $$\widehat{\mathcal{P}}(t)=(\widehat{p}_1(t):\widehat{p}_2(t):\widehat{p}(t))=$$$$=(p_1(t):p_2(t)+Z_1p_{3}(t)+\cdots +Z_{n-2}p_n(t):p(t))\in\mathbb{P}^2((\mathbb{K}(Z))(t)).$$
This notation is used for the sake of simplicity, but we note that $\widehat{\mathcal{P}}(t)$ depends on $Z$. Observe that $\widehat{\mathcal{P}}(t)$ is a proper parametrization of a rational plane curve $\widehat{{\cal C}}$ defined over the algebraic closure of  ${\Bbb K}(Z)$.

\para

We can establish a correspondence between the points of $\mathcal{C}$ and the points of $\widehat{{\cal C}}$. More precisely, for each point $P=(a_1:a_2: a_3:\cdots:a_n: a_{n+1})\in {\cal C}$ we have another point $\widehat{P}=(a_1:a_2+Z_1a_{3}+\cdots +Z_{n-2}a_n: a_{n+1}) \in \widehat{{\cal C}}$. Moreover, this correspondence is bijective for the points satisfying that $a_1\not=0$ or $a_{n+1}\not=0$. For these points, it holds that  ${\mathcal F}_{\mathcal{P}}(P)={\mathcal F}_{\widehat{\mathcal{P}}}(\widehat{P})$, which implies that $H_P(s)=H_{\widehat{P}}(s)$. Note that  the polynomial $H_{P}$  represents the fibre function of a point $P$ in the space curve $\cal C$ computed from $\mathcal{P}(t)$; i.e. the roots of $H_{P}$ are the fibre of    $P\in {\cal C}$ (this notion was introduced in Definition \ref{D-fibrefunction} for a given plane curve but it can be easily generalized for space curves).

 \para

 Observe that the above correspondence may also be established between the places of $\cal C$ and $\widehat{\cal C}$ centered at $\widehat{P}$ and $P$, respectively. That is, for each place $\varphi(t)=(\varphi_1(t):\varphi_2(t): \varphi_3(t):\cdots:\varphi_n(t): \varphi_{n+1}(t))$ of ${\cal C}$ centered at $P$ we have the place $\widehat{\varphi}(t)=(\varphi_1(t):\varphi_2(t)+Z_1\varphi_3(t)+\cdots+Z_{n-2}\varphi_n(t): \varphi_{n+1}(t))$ of $\widehat{{\cal C}}$  centered at $\widehat{P}$. Hence, the number of tangents of $\cal C$ at $P$ is the same that the number of tangents of $\widehat{\cal C}$ at $\widehat{P}$ and, as a consequence, $\mult_{{P}}({\mathcal{C}})=\mult_{\widehat{P}}(\widehat{\mathcal{C}})$.

\para

The correspondence above  introduced is not bijective if $a_1=a_{n+1}=0$. In this case, we have that  $\widehat{P}=(0:1: 0) \in \widehat{{\cal C}}$ and the corresponding points in $\cal C$ are all the points of the form $(0:a_2: a_3:\cdots:a_n: 0)$. Thus,  if we have exactly $\ell$ points $P_1,\ldots,P_\ell\in {\cal C}$ with $P_i=(0:a_{2,i}: a_{3,i}:\cdots:a_{n,i}: 0),\,i=1,\ldots,\ell$, then  ${\mathcal F}_{\widehat{\mathcal{P}}}(\widehat{P})=\cup_{i=1}^\ell {\mathcal F}_{\mathcal P}(P_i)$. Hence, $H_{\widehat{P}}(s)=\prod_{i=1}^{\ell} H_{P_i}(s)$ (note that ${\mathcal F}_{\mathcal P}(P_i)\cap{\mathcal F}_{\mathcal P}(P_j)=\emptyset$ if $i\neq j$) and $\mult_{\widehat{P}}(\widehat{\mathcal{C}})=\sum_{i=1}^\ell \mult_{P_i}(\mathcal{C})$.

\para

Thus, in order to study the singularities of $\mathcal{C}$ through those of $\widehat{\mathcal{C}}$, an additional difficulty arises when $\mathcal{C}$ contains two or more points of the form $(0:a_2: a_3:\cdots:a_n: 0)$. Let us call them {\it bad points}. In the following, we  may assume w.l.o.g. that we are not in this case, i.e.   $\mathcal{C}$ does not have two or more bad points. Otherwise, we  apply a change of coordinates, and we consider the new   parametrization
$\mathcal{P}^*(t)=(p^*_1(t):p_2(t):\cdots: p_n(t):p(t))$ of the transformed curve  $\mathcal{C}^*$, where $p^*_1=\sum_{i=1}^n\lambda_ip_i,\,\,\lambda_i\in {\Bbb K}$. By appropriately choosing   $\lambda_1,\ldots,\lambda_n\in {\Bbb K}$, we have that $\gcd(p^*_1,p)=1$ (note that $\gcd(p_1,\ldots,p_n,p)=1$) and thus,  $\mathcal{C}^*$ does not have bad points.

\para

Finally, we also note that additional points, which can not be written in the form $(a_1:a_2+Z_1a_{3}+\cdots +Z_{n-2}a_n: a_{n+1}), \,a_i\in {\Bbb K},\,i=1,\ldots,n+1$, may appear in the curve  $\widehat{{\cal C}}$. Such points are obtained as $\widehat{\mathcal{P}}(t)$ for $t\in\mathbb{K}(Z)\setminus \mathbb{K}$ and they do not have a correspondence with any point of $\cal C$.

\para

Under these conditions, let $\widehat{G}_1$, $\widehat{G}_2$ and $\widehat{G}_3$ be the equivalent polynomials to $G_1$, $G_2$ and $G_3$ (defined in (\ref{Eq-fibra-generica})), constructed from the parametrization $\widehat{\mathcal{P}}(t)$. In addition, let $\widehat{\delta}_i:=\degree_t(\widehat{G}_i)$ and $\widehat{\lambda}_{ij}:=\min\{\widehat{\delta}_i,\widehat{\delta}_j\},\,i,j=1,2,3,\,i<j$, $$\widehat{G}_i^*(s,t):=\displaystyle\frac{\widehat{G}_i(s,t)}{t-s}\in ({\Bbb K}[Z])[s,t],\,\,\,i=1,2,3,$$
and
$$\widehat{R}_{ij}(s):=\Res_t(\widehat{G}_i^*,\widehat{G}_j^*)\in ({\Bbb K}[Z])[s],\,\,\, i,j=1,2,3,\,i<j.$$

\para

\noindent
Then, the T--function of the parametrization  $\widehat{\mathcal{P}}(t)$
is given by
$$\widehat{T}(s)=\widehat{R}_{12}(s)/\widehat{p}(s)^{\widehat{\lambda}_{12}-1}.$$
It holds that $\widehat{T}(s)\in ({\Bbb K}[Z])[s]$  (see Corollary \ref{C-T-polin}), and by Proposition \ref{P-res13-res23} we get that
$$\widehat{T}(s)=\frac{\widehat{R}_{12}(s)}{\widehat{p}(s)^{\widehat{\lambda}_{12}-1}}=\frac{\widehat{R}_{13}(s)}{\widehat{p}_1(s)^{\widehat{\lambda}_{13}-1}}=\frac{\widehat{R}_{23}(s)}{\widehat{p}_2(s)^{\widehat{\lambda}_{23}-1}}.$$

\para

The following theorem is obtained as a consequence of Theorem \ref{T-tfunct} (see Section \ref{S-resultantsingularities}), and it shows how the function $\widehat{T}(s)$ can be used to define an equivalent polynomial to the T--function  introduced for plane curves (see Definition \ref{D-T-funct}). This polynomial, will be denoted as $T_E(s)$.
%,  describes totally the singularities of $\cal C$,   since the factorization of  $T_E(s)$  provides the fibre functions of each singularity as well as its corresponding multiplicity.

\para

Similarly as in the case of  plane curves, we assume that the space curve,   $\cal C$,  has only ordinary singularities.  The case of space curves with non--ordinary singularities will be analyzed in a future work and an equivalent theorem to Theorem \ref{T-tfunct-spa} will be obtained for this case.

\para

Finally we remind, that if  $\mathcal{C}$  has two or more bad points, we consider the transformed curve $\mathcal{C}^*$ introduced above. Note that $H_P(s)=H_{P^*}(s)$, where $P\in {\cal C}$ is moved to  the point $P^*\in \mathcal{C}^*$ when the change of coordinates is applied. Undoing this change of coordinates, one recovers the initial singularities $P\in {\cal C}$.

\para

\begin{theorem}\label{T-tfunct-spa}
Let $\mathcal{C}$ be a rational algebraic space curve defined by a parametrization $\mathcal{P}(t)$, with limit point $P_L$. Let $P_1,\ldots,P_n$ be the singular points of $\cal C$, with multiplicities $m_1,\ldots,m_n$ respectively. Let us assume that they are ordinary singularities and that $P_i\neq P_L$, for $i=1,\ldots,n$. Then, it holds that
$$T_E(s)=\prod_{i=1}^nH_{P_i}(s)^{m_i-1},$$
where $T_E(s)=\Content_Z\left(\widehat{T}(s)\right)\in {\Bbb K}[s],$ and $\Content_Z(\widehat{T})$ represents the content of the polynomial $\widehat{T}$ w.r.t $Z$.
\end{theorem}

\noindent\textbf{Proof:} From the above statements, we observe that there exists a bijective correspondence between the points $\widehat{P}=(a_1:a_2+Z_1a_{3}+\cdots +Z_{n-2}a_n: a_{n+1}), \,a_i\in {\Bbb K},\,i=1,\ldots,n+1,$  of $\widehat{{\cal C}}$ and the points  $P=(a_1:a_2: a_3:\cdots:a_n: a_{n+1})$  of ${\cal C}$. Consequently, we have that $\mult_{\widehat{P}}(\widehat{\mathcal{C}})=\mult_{{P}}({\mathcal{C}})$, which implies that $\widehat{P}$ is a singularity of $\widehat{{\cal C}}$ of multiplicity $m$ if and only if $P$ is a singularity of ${\cal C}$ of multiplicity $m$. Hence, using  Theorem \ref{T-tfunct}, we deduce that \[\widehat{T}(s)=\prod_{i=1}^nH_{P_i}(s)^{m_i-1}L(s,Z).\]
We observe that the factor $L(s,Z)\in {\Bbb K}[s,Z]\setminus{\Bbb K}[s]$ is a product of  the fibre functions corresponding to the singularities of  $\widehat{{\cal C}}$ that can not be written as $(a_1:a_2+Z_1a_{3}+\cdots +Z_{n-2}a_n: a_{n+1}), \,a_i\in {\Bbb K},\,i=1,\ldots,n+1$ (these singularities do not have an equivalent singularity in $\cal C$, and its fibre function necessarily is a polynomial in ${\Bbb K}[s,Z]\setminus{\Bbb K}[s]$). Then, we  conclude that
  $$T_E(s)=\Content_Z\left(\widehat{T}(s)\right)=\prod_{i=1}^nH_{P_i}(s)^{m_i-1}.$$\hfill $\Box$
%In order to discard these singularities, we have to consider the content w.r.t $Z$ in the polynomial $\widehat{T}(s)$

\begin{example}
Let $\mathcal{C}$ be the rational space curve defined by the projective parametrization
$\mathcal{P}(t)=(p_1(t):p_2(t):p_3(t):p(t))\in\mathbb{P}^3(\mathbb{C}(t)),$
where
$$\begin{array}{l}
    p_1(t)=t^5-5t^4+5t^3+5t^2-6t \\
    p_2(t)=t^5+5t^4+5t^3-5t^2-6t \\
    p_3(t)=t^5+7t^4+17t^3+17t^2+6t\\
    p(t)=t^4-13t^2+36.
  \end{array}$$
We consider the plane curve $$\widehat{\mathcal{P}}(t)=(p_1(t):p_2(t)+Zp_{3}(t):p(t))\in\mathbb{P}^2((\mathbb{C}(Z))(t))$$
and compute the corresponding T--function:
$$\widehat{T}(s)=298598400(s-2)(s-3)(s+3)(s+2)s(s+1)L(s,Z),$$
where $L(s,Z)=(3Z+1)(2Z+1)(25s^6+25Z^2s^6+50Zs^6+60Z^2s^5+35Zs^5-25s^5+125s^4+322Z^2s^4+375Zs^4+360Z^2s^3-65Zs^3-125s^3-150s^2-185Zs^2+229Z^2s^2+150Zs+300Z^2s+150s-360Z+432Z^2)$. Removing $L(s,Z)$ (which depends on $Z$), we get
$$T_E(s)=298598400(s-2)(s-3)(s+3)(s+2)s(s+1).$$

\begin{figure}[h]
$$
\psfig{figure=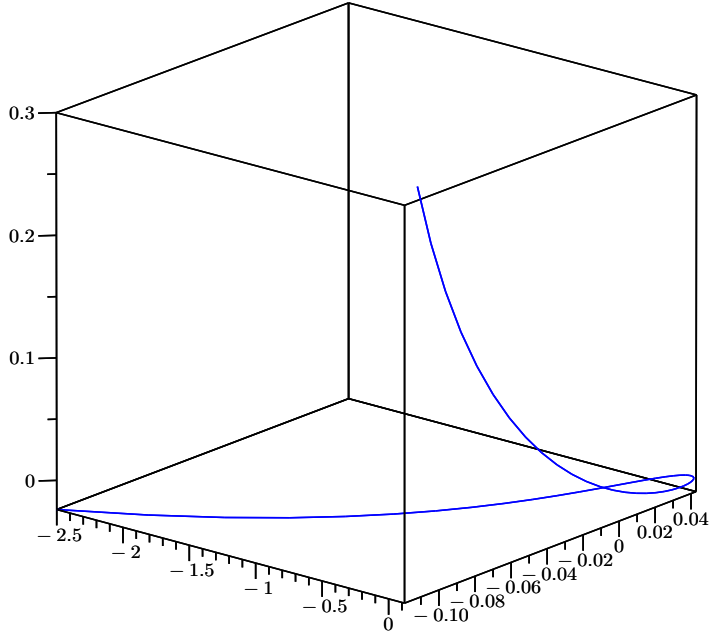,width=12cm,height=10cm}
$$ \caption{Curve $\cal C$ and the double point $P_3$}\label{F-espacial}
\end{figure}
\noindent Now, reasoning as in Example \ref{E-Tfunct-planas}, we deduce that $\cal C$ has three singularities:

\begin{itemize}
\item The infinity point $P_1=(0:1:3:0)$, with fibre function $H_{P_1}(t)=(t-2)(t-3)$ (let us remark that $P_1$ is a bad point; however, it does not represent a problem in this case since there are no more bad points in $\mathcal{C}$).%(note that, in this case, $\widehat{P}_1=(0:1:0)$ and $P_1=\cP(2)=\cP(3)$; thus $H_{\widehat{P}_1}(t)=H_{P_1}(t)$).
\item The infinity point $P_2=(1:0:0:0)$, with fibre function $H_{P_2}(t)=(t+2)(t+3).$
\item The affine point $P_3=(0:0:0:1)$, with fibre function $H_{P_3}(t)=t(t+1).$
\end{itemize}
Note that $P_1,\,P_2$ and $P_3$ are double points of $\cal C$ (see Figure \ref{F-espacial}). %One may check that $\cal C$ only has an affine singularity in $(0,0,0)$.
\end{example}

\section{Proof of the main theorem}\label{S-proof}

This section is devoted to show the main result of this paper, Theorem \ref{T-tfunct} in Section \ref{S-resultantsingularities}. For this purpose, we first prove some previous results. In particular, the following lemma is obtained using the main properties of the resultants (see e.g. \cite{Cox1998}, \cite{SWP}, \cite{Vander}).

\para

\begin{lemma}\label{L-props-resultante}
Let $A(s,t), B(s,t), C(s,t)\in {\Bbb K}[s,t]$, and $K(s)\in {\Bbb K}[s]$. The following properties hold: %{\color{red}Up to constants in ${\Bbb K}\setminus \{0\}$,}
\begin{enumerate}
\item $\Res_t(A,K)=K^{\deg_t(A)}$.
\item $\Res_t(A,B\cdot C)=\Res_t(A,B)\cdot \Res_t(A,C).$
\item If $B$ divides  $A$, it holds that $\Res_t(A/B,C)=\Res_t(A,C)/\Res_t(B,C).$
\item $\Res_t(A,B+CA)=\lc(A)^k\Res_t(A,B)$, where $k={\deg_t(B+CA)-\deg_t{B}}$.
%\item Sea $M(t,s)=A(t,s)C(t)+B(s)D(t)$. Entonces
%$$Res_t(M(t,s),C(t))=lc_t(C(t))^{\deg_t(M)-\deg(D)}\cdot B(s)^{\deg(C)}\cdot Res_t(C(t),D(t)).$$
\end{enumerate}
\end{lemma}

\noindent\textbf{Proof:} First, we remind that if $A(t), B(t)\in {\Bbb K}[t]$, it holds that
$$\Res_t(A,B)=\lc(A)^{\deg(B)}\prod_{A(\alpha_i)=0}B(\alpha_i).$$
Now, let  $A(s,t)\in ({\Bbb K}[s])[t]$. The leader coefficient of $A(s,t)$ w.r.t the variable $t$ is $\lc_t(A)\in{\Bbb K}[s]$ and, for each $s\in\mathbb{K}$,  the polynomial $A(s,t)$ has $\deg_t(A)$ roots, $\alpha_1(s),\ldots,\alpha_{\deg_t(A)}(s)$, in the algebraic closure of ${\Bbb K}(s)$, such that $A(s,\alpha_i(s))=0$, for $i=1,\ldots,\deg_t(A)$. Then,
\begin{equation}\label{Eq-resultant}\Res_t(A,B)=\lc_t(A)^{\deg_t(B)}\prod_{A(s,\alpha_i(s))=0}B(s,\alpha_i(s)).\end{equation}
Now, we prove  the four statements of the lemma.
\begin{enumerate}
\item The first statement follows using (\ref{Eq-resultant}) for the case that $B(s,t)\in{\Bbb K}[s]$. Then $\deg_t(B)=0$ and $B(s,\alpha_i(s))=B(s)$ for each $i=1,\ldots,\deg_t(A)$.
\item In order to prove statement $2$, we use (\ref{Eq-resultant}), and we get that
$$\Res_t(A,B\cdot C)=\lc_t(A)^{\deg_t(B\cdot C)}\prod_{A(s,\alpha_i(s))=0}B(s,\alpha_i(s))\cdot C(s,\alpha_i(s)).$$
Since $\deg_t(B\cdot C)=\deg_t(B)+\deg_t(C)$, we have that
$$\Res_t(A,B\cdot C)=$$$$\left(\lc_t(A)^{\deg_t(B)}\prod_{A(s,\alpha_i(s))=0}B(s,\alpha_i(s))\right)\left(\lc_t(A)^{\deg_t(C)}\prod_{A(s,\alpha_i(s))=0}C(s,\alpha_i(s))\right)$$$$=\Res_t(A,B)\cdot \Res_t(A,C).$$
\item Reasoning similarly as in statement $2$, we get that
$$\Res_t(A/B,C)=\lc_t(A/B)^{\deg_t(C)}\prod_{A(s,\alpha_i(s))=0,\,B(s,\alpha_i(s))\neq 0}C(s,\alpha_i(s)).$$
Since $B$ divides  $A$, we obtain that
$$\lc_t(A/B)^{\deg_t(C)}\prod_{A(s,\alpha_i(s))=0,\,B(s,\alpha_i(s))\neq 0}C(s,\alpha_i(s))=$$$$\frac{\lc_t(A)^{\deg_t(C)}}{\lc_t(B)^{\deg_t(C)}}\frac{\prod_{A(s,\alpha_i(s))=0}C(s,\alpha_i(s))}{\prod_{B(s,\alpha_i(s))=0}C(s,\alpha_i(s))}=\frac{\Res_t(A,C)}{\Res_t(B,C)}.$$
\item We reason similarly as above, and  we get that
$$\Res_t(A,B+CA)=\lc_t(A)^{\deg_t(B+CA)}\cdot$$$$\prod_{A(s,\alpha_i(s))=0}(B(s,\alpha_i(s))+C(s,\alpha_i(s))A(s,\alpha_i(s)))=$$
$$=\lc(A)^{\deg_t(B+CA)-\deg_t{B}}\left(\lc(A)^{\deg_t(B)}\prod_{A(s,\alpha_i(s))=0}B(s,\alpha_i(s))\right)=$$$$=\lc(A)^k\Res_t(A,B),\,\,\mbox{ where } k={\deg_t(B+CA)-\deg_t{B}}.$$
%\item Recordemos que $Res_t(M(t,s),C(t))=Res_t(C(t),M(t,s))$ salvo constantes. Entonces,
%$$Res_t(C(t),M(t,s))=lc_t(C(t))^{\deg_t(M)}\prod_{C(\alpha_i)=0}(A(\alpha_i,s)C(\alpha_i)+B(s)D(\alpha_i))$$
%$$=lc_t(C(t))^{\deg_t(M)}\prod_{C(\alpha_i)=0}B(s)D(\alpha_i)=lc_t(C(t))^{\deg_t(M)}B(s)^{\deg(C)}\prod_{C(\alpha_i)=0}D(\alpha_i)$$
%$$=lc_t(C(t))^{\deg_t(M)-\deg_t(D)}B(s)^{\deg(C)}\left(lc_t(C(t))^{\deg_t(D)}\prod_{C(\alpha_i)=0}D(\alpha_i)\right).$$
%Pero ahora lo que hay entre parentesis coincide ser $Res_t(C(t),D(t))$, de donde se deduce el resultado.\hfill $\Box$
\end{enumerate}

\para

The following lemma provides a first approach to the proof of the main result presented in this paper (Theorem \ref{T-tfunct} in Section \ref{S-resultantsingularities}).  In particular, it is shown that each factor $H_P(s)^{m-1}$, where $P$ is an ordinary singular point of multiplicity $m$, divides the T--function.

\para

\begin{lemma}\label{L-tfunct} Let $\mathcal{C}$ be a rational algebraic curve defined by a parametrization $\mathcal{P}(t)$, with limit point $P_L$. Let $P\neq P_L$  be an ordinary singular point of multiplicity $m$. It holds that
$$ T(s)=H_P(s)^{m-1}T^*(s),$$
where $T^*(s)\in {\Bbb K}[s]$ and $\gcd(H_P(s),T^*(s))=1.$
\end{lemma}
\noindent\textbf{Proof:} In order to prove this lemma, we distinguish three different steps. In  Step 1, we prove that the lemma holds if   $P=(0:0:1)$. In Step 2, we show that  the lemma holds for any affine singularity. Finally, in Step 3, we prove the lemma for $P$ being a singular point at infinity.

\begin{center}
{\it Step 1}
\end{center}
Let us assume that the given singularity is the point
$P=(0:0:1)$. Note that
$H_P(t)=\gcd(\phi_1,\phi_2)=\gcd(p_1,p_2)$ since,  in this case,
$a=b=0$ (see (\ref{Eq-fibre-equations})). Thus, we may write
$$\left\{\begin{array}{l}p_1(t)=H_P(t)\overline{p}_1(t)\\p_2(t)=H_P(t)\overline{p}_2(t),\end{array}\right.$$
where $\overline{p}_1(t)$ and $\overline{p}_2(t)$ are polynomials satisfying that
$\gcd(\overline{p}_1,\overline{p}_2)=1$. In addition, it holds that  $\gcd(H_P(t),p(t))=1$, since   $\gcd(p_1,p_2,p)=1$. Hence, from (\ref{Eq-fibra-generica}), we may write
$G_3(s,t)=H_P(s)H_P(t)(\overline{p}_1(s)\overline{p}_2(t)-\overline{p}_2(s)\overline{p}_1(t))$ that is,
\begin{equation}\label{Eq-aux-G3}G_3(s,t)=H_P(s)H_P(t)\overline{G}_3(s,t),\end{equation} where
$\overline{G}_3(s,t)=\overline{p}_1(s)\overline{p}_2(t)-\overline{p}_2(s)\overline{p}_1(t)$.

\para

Observe that since $P\neq P_L$, Corollary
\ref{C-multiplicidad} holds and then  $\deg(H_P(t))=m\geq 2$. Hence,  there exist at least two values $s_0, s_1\in {\Bbb K}$  such that $H_P(s_0)=H_P(s_1)=0$  and then, since these roots belong to the fibre of  $P$, we deduce that
$$\left(\left(\frac{p_1}{p}\right)(s_0),\left(\frac{p_2}{p}\right)(s_0)\right)=\left(\left(\frac{p_1}{p}\right)(s_1),\left(\frac{p_2}{p}\right)(s_1)\right)=(0,0).$$
Since $P$ is an ordinary singularity, we have that there can not exist $K_1, K_2\in {\Bbb K}$ such that
\begin{equation}\label{Eq-aux-noprop}K_1\left(\left(\frac{p_1}{p}\right)'(s_0),\left(\frac{p_2}{p}\right)'(s_0)\right)=K_2\left(\left(\frac{p_1}{p}\right)'(s_1),\left(\frac{p_2}{p}\right)'(s_1)\right).\end{equation}

\para

\noindent
We also note that, for $i=1,2$ and $j=0,1$, it holds that
$$\left(\frac{p_i}{p}\right)'(s_j)=\frac{p'_i(s_j)p(s_j)-p_i(s_j)p'(s_j)}{p(s_j)^2}=\frac{p'_i(s_j)}{p(s_j)},$$
(note that $p_i(s_j)=0$).  In addition, since  $p_i(t)=H_P(t)\overline{p}_i(t)$, we get that
$$\frac{p'_i(s_j)}{p(s_j)}=\frac{H'_P(s_j)\overline{p}_i(s_j)+H_P(s_j)\overline{p}'_i(s_j)}{p(s_j)}=\frac{H'_P(s_j)\overline{p}_i(s_j)}{p(s_j)},$$
(note that $H_P(s_j)=0$). By substituting in  (\ref{Eq-aux-noprop}), we obtain that
$$K_1\left(\frac{H'_P(s_0)\overline{p}_1(s_0)}{p(s_0)},\frac{H'_P(s_0)\overline{p}_2(s_0)}{p(s_0)}\right)\neq K_2\left(\frac{H'_P(s_1)\overline{p}_1(s_1)}{p(s_1)},\frac{H'_P(s_1)\overline{p}_2(s_1)}{p(s_1)}\right).$$
Hence, we remark that:
\begin{itemize}
\item $H'_P(s_i)\neq 0,\,i=0,1$. That is,  $s_0$ and $s_1$ are simple roots of  $H_P$.
\item None of the following equalities may be verified:
$$\left\{\begin{array}{l}\overline{p}_1(s_0)=\overline{p}_2(s_0)=0\\
\overline{p}_1(s_1)=\overline{p}_2(s_1)=0\\
\overline{p}_1(s_0)=\overline{p}_1(s_1)=0\\
\overline{p}_2(s_0)=\overline{p}_2(s_1)=0\end{array}\right.$$
%Neither $\overline{p}_1(s_i)=\overline{p}_2(s_i)=0, i=0,1,$ nor $\overline{p}_i(s_0)=\overline{p}_i(s_1)=0, i=1,2,$ may be verified.
% and $\overline{p}_1(s_1)=\overline{p}_2(s_1)=0$ and $\overline{p}_1(s_0)=\overline{p}_1(s_1)=0$ and $\overline{p}_2(s_0)=\overline{p}_2(s_1)=0$.
\item If $\overline{p}_1(s_1)\neq 0$ and $\overline{p}_2(s_1)\neq 0$, then \begin{equation}\label{Eq-aux-noord}\frac{\overline{p}_1(s_0)}{\overline{p}_1(s_1)}\neq\frac{\overline{p}_2(s_0)}{\overline{p}_2(s_1)}.\end{equation}
\end{itemize}

Taking into account these remarks, we prove that $T(s)=H_P(s)^{m-1}T^*(s)$. For this purpose, we first write the T--function as
$$T(s)=R_{13}(s)/p_1(s)^{\lambda_{13}-1} $$
(see Proposition \ref{P-res13-res23}). From (\ref{Eq-aux-G3}), we get that
$$R_{13}(s)=\Res_t\left(G^*_1(s,t),H_P(s)H_P(t)\overline{G}_3^*(s,t)\right),$$
where $\overline{G}_3^*(s,t)=\overline{G}_3(s,t)/(t-s)$ (note that $\overline{G}_3^*(s,t)\in {\Bbb K}[s,t]$ since $(t-s)$ divides
$\overline{G}_3(s,t)$). By applying statement 2 of Lemma
\ref{L-props-resultante}, we have that $$R_{13}(s)=$$
\begin{equation}\label{Eq-aux-R}\Res_t\left(G^*_1(s,t),H_P(s)\right)\Res_t\left(G^*_1(s,t),H_P(t)\right)\Res_t\left(G^*_1(s,t),\overline{G}_3^*(s,t)\right).\end{equation}
Let us analyse the first two factors:
\begin{itemize}
\item From statement 1 of Lemma \ref{L-props-resultante}, we have that
$$\Res_t\left(G^*_1(s,t),H_P(s)\right)=H_P(s)^{\deg_t(G^*_1)}=H_P(s)^{\delta_1-1}.$$
\item On the other side,
$$\Res_t\left(G^*_1(s,t),H_P(t)\right)=\Res_t\left(\frac{{p}_1(s){p}(t)-{p}(s){p}_1(t)}{t-s},H_P(t)\right)$$
and, by applying statement 3 of Lemma  \ref{L-props-resultante}, we get that
$$\frac{\Res_t\left({p}_1(s){p}(t)-{p}(s){p}_1(t),H_P(t)\right)}{\Res_t\left({t-s},H_P(t)\right)}.$$
Note that $\Res_t\left({t-s},H_P(t)\right)=H_P(s)$ and, since  $p_1(t)=H_P(t)\overline{p}_1(t)$, the above expression can be written as
$$\frac{\Res_t\left({p}_1(s){p}(t)-{p}(s)H_P(t)\overline{p}_1(t),H_P(t)\right)}{H_P(s)}.$$
Now, from statements 1 and 4 of Lemma  \ref{L-props-resultante}, we get that
$$\Res_t\left(G^*_1(s,t),H_P(t)\right)=\frac{p_1(s)^m \Res_t\left(p(t),H_P(t)\right)\lc(H_P(t))^k}{H_P(s)},$$
where $k\in\mathbb{K}$. Note that,  $\lc(H_P(t))^k\in{\Bbb K}\setminus\{0\}$ and   $\Res_t\left(p(t),H_P(t)\right)\in{\Bbb K}\setminus\{0\}$ (since $\gcd(p,H_P)=1$). Furthermore, we have that
$p_1(s)=H_P(s)\overline{p}_1(s)$. Thus, up to constants in  ${\Bbb K}\setminus\{0\}$, we deduce that
$$\Res_t\left(G^*_1(s,t),H_P(t)\right)=H_P(s)^{m-1}\overline{p}_1(s)^m.$$
\end{itemize}
Substituting in (\ref{Eq-aux-R}), we obtain that
$$R_{13}(s)=H_P(s)^{\delta_1-1}H_P(s)^{m-1}\overline{p}_1(s)^m \Res_t\left(G^*_1(s,t),\overline{G}_3^*(s,t)\right),$$
and thus
$$T(s)=\frac{R_{13}(s)}{p_1(s)^{\lambda_{13}-1}}=\frac{R_{13}(s)}{H_P(s)^{\lambda_{13}-1}\overline{p}_1(s)^{\lambda_{13}-1}}=H_P(s)^{m-1}T^*(s),$$
where
$$T^*(s)=\frac{\Res_t\left(G^*_1(s,t),\overline{G}_3^*(s,t)\right)H_P(s)^{\delta_1-\lambda_{13}}}{\overline{p}_1(s)^{\lambda_{13}-1-m}}.$$
Note that $\lambda_{13}=\delta_1$ since $\delta_1\leq\delta_3$. Otherwise, if $\delta_1>\delta_3$, we would have that  $\max\{d_1,d_3\}>\max\{d_1,d_2\}$ and then, $d_3>d_1,d_2$. However, this would imply that  $P=P_L$ (see  Definition \ref{D-pto-limite}), which contradicts the assumption. Therefore,
\begin{equation}\label{Eq-T-ast1}
T^*(s)=\frac{\Res_t\left(G^*_1(s,t),\overline{G}_3^*(s,t)\right)}{\overline{p}_1(s)^{\delta_1-1-m}}.
\end{equation}
Now, we  prove that   $T^*(s) \in {\Bbb K}[s]$. We can assume that $\delta_1-1-m\geq 0$, since $\delta_1=\max\{d_1,d_3\}\geq d_1$ and $d_1=\deg(p_1)=\deg(H_P\cdot\overline{p}_1)=m+\deg(\overline{p}_1)$. Hence,  $\delta_1\geq m+1$ except for the case that $\deg(\overline{p}_1)=0$, but in this situation we would have that
$$T^*(s)= {\Res_t\left(G^*_1(s,t),\overline{G}_3^*(s,t)\right)}\in {\Bbb K}[s].$$
So, let $\delta_1-1-m\geq 0$. Now, we reason similarly as in the proof of Corollary  \ref{C-T-polin}. Indeed: let us assume that  $T^*(s)$ is not a polynomial. Then, by taking the simplified rational function, we may write
$$T^*(s)=\frac{M_{13}(s)}{\widehat{p}_1(s)}$$
where $M_{13}(s)\in {\Bbb K}[s]$,\,\,$\widehat{p}_1(s)\in {\Bbb K}[s]\setminus {\Bbb K}$  and $\gcd(M_{13}(s),\widehat{p}_1(s))=1$. Note that $\widehat{p}_1(s)$ divides $\overline{p}_1(s)^{\delta_1-1-m}$.
%, luego los distintos factores de $\widehat{p}_1(s)$ deben dividir a $\overline{p}_1(s)$.

\para

We observe that  (\ref{Eq-T-ast1}) is obtained from
$T(s)=R_{13}(s)/p_1(s)^{\lambda_{13}-1}.$
However, taking into account Proposition \ref{P-res13-res23}, we could have considered the expression
$$T(s)=R_{23}(s)/p_2(s)^{\lambda_{23}-1}$$
concluding that  $T(s)=H_P(s)^{m-1}T^*(s)$,
where
\begin{equation}\label{Eq-T-ast2}
T^*(s)=\frac{\Res_t\left(G^*_2(s,t),\overline{G}_3^*(s,t)\right)}{\overline{p}_2(s)^{\delta_2-1-m}}.
\end{equation}
Reasoning similarly as above, we get that there would exist $M_{23}(s)\in {\Bbb K}[s]$,\,$\widehat{p}_2(s)\in {\Bbb K}[s]\setminus {\Bbb K}$  with  $\gcd(M_{23}(s),\widehat{p}_2(s))=1$, such that
$$T^*(s)=\frac{M_{23}(s)}{\widehat{p}_2(s)}.$$
In addition,  $\widehat{p}_2(s)$ would divide $\overline{p}_2(s)^{\delta_2-1-m}$. Thus, we have that
$$\frac{M_{13}(s)}{\widehat{p}_1(s)}=\frac{M_{23}(s)}{\widehat{p}_2(s)}$$
where  $\gcd(M_{13}(s),\widehat{p}_1(s))=\gcd(M_{23}(s),\widehat{p}_2(s))=1$, which implies that  $\widehat{p}_1(s)=\widehat{p}_2(s)$. Hence $\gcd(\overline{p}_1,\overline{p}_2)\neq 1$, which contradicts the definition of these functions.

\para

Thus, we have proved that $T^*$ is a polynomial. Finally, we show that  $\gcd(H_P(s),T^*(s))=1$. Indeed: if $\gcd(H_P(s), T^*(s))\not=1$, there exists  $s_0\in {\Bbb K}$ such that $H_P(s_0)=0$ and
$T^*(s_0)=0$, which implies that $$\Res_t\left(G^*_1(s,t),\overline{G}_3^*(s,t)\right)(s_0)=0.$$
Then, by the properties of the resultants, one of the following statements hold:
\begin{enumerate}
\item  There exists  $s_1\in {\Bbb K}$ such that $G^*_1(s_0,s_1)=\overline{G}_3^*(s_0,s_1)=0$. This would imply that     $G^*_3(s_0,s_1)=0$, and, then, $s_0$ and $s_1$ are elements of the fibre of  $P$. On the other side, $$\overline{G}_3(s_0,s_1)=\overline{p}_1(s_0)\overline{p}_2(s_1)-\overline{p}_2(s_0)\overline{p}_1(s_1)=0$$ and thus,
$$\frac{\overline{p}_1(s_0)}{\overline{p}_1(s_1)}=\frac{\overline{p}_2(s_0)}{\overline{p}_2(s_1)}.$$
This implies that  $P$ is a non--ordinary singular point
(see (\ref{Eq-aux-noord})), which contradicts the assumptions.
\item It holds that $\gcd(\lc_t(G^*_1), \lc_t(\overline{G}_3^*))(s_0)=0$. Then, in particular,
$$\lc_t(G^*_1)(s_0)=\lc_t(G_1)(s_0)=p_1(s_0)c_d-p(s_0)a_d=0\Rightarrow a_d=0.$$
Now we reason similarly with the equality  $T(s)=R_{23}(s)/p_2(s)^{\lambda_{23}-1}$ and we get  (\ref{Eq-T-ast2}). From this expression, and reasoning similarly as above, we obtain that  $\gcd(\lc_t(G^*_2), \lc_t(\overline{G}_3^*))(s_0)=0$, which implies that $b_d=0$. However, if $a_d=b_d=0$ we deduce that $P=P_L$, which contradicts the assumptions.

%\item It holds that $\gcd(\lc_t(G^*_1), \lc_t(\overline{G}_3^*))(s_0)=0$. In this case,
%$$\lc_t(G^*_1)(s_0)=\lc_t(G_1)(s_0)=p_1(s_0)c_d-p(s_0)a_d=0\Rightarrow a_d=0,$$
%and
%$$\lc_t(\overline{G}^*_3)(s_0)=\lc_t(\overline{G}_3)(s_0)=\overline{p}_1(s_0)b_d-\overline{p}_2(s_0)a_d=0\Rightarrow
%\overline{p}_1(s_0)=0.$$
%Under these conditions, we may reason similarly with the equality  $T(s)=R_{23}(s)/p_2(s)^{\lambda_{23}-1}$, and we would get  (\ref{Eq-T-ast2}). From this expression, and reasoning similarly as above, we obtain that  $\gcd(\lc_t(G^*_2), \lc_t(\overline{G}_3^*))(s_0)=0$, which implies that $b_d=0$ and $\overline{p}_2(s_0)=0$. However, if $a_d=b_d=0$ we deduce that $P=P_L$, which contradicts the assumptions.
\end{enumerate}
Therefore, we conclude that  $\gcd(H_P(s),T^*(s))=1$.

\begin{center}
{\it Step 2}
\end{center}

Let $P=(a:b:1)$ be a singularity of multiplicity  $m$. In this case, we consider the translation of the curve  $\cal C$ defined by the parametrization
$$\widetilde{\cal P}(t)=(p_1(t)-ap(t):p_2(t)-bp(t):p(t)).$$
We have that the point  $P=(a:b:1)$ moves to the point $\widetilde{P}=(0:0:1)$, and then
$H_P(t)=H_{\widetilde{P}}(t)$ (note that the polynomial  $H_P(t)$ is computed from ${\cal P}(t)$, and the polynomial $H_{\widetilde{P}}(t)$ is computed from $\widetilde{\cal P}(t)$).
%(aqui
%estamos cometiendo un cierto abuso de notacion, ya que una y otra
%funcion de fibra se obtienen a partir de parametrizaciones
%distintas).

\para

On the other side, if we compute the polynomial equivalent to  $G_1(s,t)$ with the new parametrization $\widetilde{\cal P}(t)$, we get that
$$\widetilde{G}_1(s,t)=\widetilde{p}_1(s)\widetilde{p}(t)-\widetilde{p}(s)\widetilde{p}_1(t)=$$
$$=(p_1(s)-ap(s))p(t)-p(s)(p_1(t)-ap(t))=p_1(s)p(t)-p(s)p_1(t)=G_1(s,t).$$
Similarly, one obtains that  $\widetilde{G}_2(s,t)=G_2(s,t)$. Thus,
$$\widetilde{R}_{12}(s)=\Res_t\left(\frac{\widetilde{G}_1(s,t)}{t-s},\frac{\widetilde{G}_2(s,t)}{t-s}\right)=\Res_t\left(\frac{G_1(s,t)}{t-s},\frac{G_2(s,t)}{t-s}\right)=R_{12}(s),$$
and then
$$\widetilde{T}(s)=\widetilde{R}_{12}(s)/\widetilde{p}(s)^{\lambda_{12}-1}=R_{12}(s)/p(s)^{\lambda_{12}-1}=T(s).$$
Thus, it holds that  $T(s)=H_P(s)^{m-1}T^*(s)$ and
$\gcd(H_P,T^*)=1$, since from Step  1, these equalities hold for  $H_{\widetilde{P}}(s)$ and $\widetilde{T}(s)$.

\begin{center}
{\it Step 3}
\end{center}

Let us prove that the lemma holds for a singularity at infinity. For this purpose, we assume that   $P=(1:0:0)$. Note that we can reason similarly as in Step 1 taking into account that for this case,
$H_P(t)=\gcd(p(t),\phi_3(t))=\gcd(p(t),p_2(t))$ (see Remark
\ref{R-afin-noafin}). Hence,
$$G_2(s,t)=H_P(s)H_P(t)(\overline{p}_2(s)\overline{p}(t)-\overline{p}(s)\overline{p}_2(t)).$$
and
$$R_{12}(s)=\Res_t\left(G^*_1,G^*_2\right)=\Res_t\left(G^*_1(s,t),H_P(s)H_P(t)\overline{G}_2^*(s,t)\right)$$
where
$$\overline{G}_2^*(s,t)=\frac{\overline{p}_2(s)\overline{p}(t)-\overline{p}(s)\overline{p}_2(t)}{t-s}.$$
Thus, using the expression $T(s)=R_{12}(s)/p(s)^{\lambda_{12}-1}$, we deduce that the lemma also holds if the singularity is the point $(1:0:0)$. A similar reasoning with the expression $T(s)=R_{23}(s)/p_2(s)^{\lambda_{23}-1}$ shows that the lemma holds for the point $(1:0:0)$.

\para

Finally, let us assume that $P=(a:b:0)$. Then, we reason similarly as in Step  2 and we apply a translation such that the point $P$ is moved to the point
$\widetilde{P}=(1:0:0)$. This translation can be defined parametrically by
$$\widetilde{\mathcal{P}}(t)=(p_1(t)):p_2(t)-(b/a)p_1(t):p(t)).$$
We assume that  $a\neq 0$; otherwise, it should be $b\neq 0$ and we would use a translation that would move   $P$ to
$(0:1:0)$.

\para

Under these conditions, we have that
$H_P(t)=H_{\widetilde{P}}(t)$.
%(igual que
%antes, estamos cometiendo un abuso de notacion, ya que $H_P$ y
%$H_{\widetilde{P}}$ se obtienen a partir de parametrizaciones
%distintas).
In addition, if we compute the equivalent polynomials to ${G}_1(s,t)$ and ${G}_3(s,t)$ with the new parametrization $\widetilde{\mathcal{P}}(t)$, we get that  $\widetilde{G}_1(s,t)=G_1(s,t)$ and
$\widetilde{G}_3(s,t)=G_3(s,t)$. Thus, from Proposition  \ref{P-res13-res23},
$$\widetilde{T}(s)=\widetilde{R}_{13}(s)/\widetilde{p}_1(s)^{\lambda_{13}-1}=R_{13}(s)/p_1(s)^{\lambda_{13}-1}=T(s).$$
Therefore $T(s)=H_P(s)^{m-1}T^*(s)$ and $\gcd(H_P,T^*)=1$, since both equalities hold for
$H_{\widetilde{P}}(s)$ and $\widetilde{T}(s)$.
 \hfill $\Box$

\para

\para

\para

 \noindent
 {\bf Proof of the Main Formula (Theorem \ref{T-tfunct} in Section \ref{S-resultantsingularities})}\\
Taking into account Lemma \ref{L-tfunct}, we have that for each singular point  $P_i$, it holds that
$T(s)=H_{P_i}(s)^{m_i-1}T_i^*(s)$, where $T_i^*(s)\in {\Bbb K}[s]$ and   $\gcd(H_{P_i},T_i^*)=1$. In addition,
$\gcd(H_{P_i},H_{P_j})=1$ for $i\neq j$ (otherwise, there would exist $s_1\in {\Bbb K}$ such that ${\cal P}(s_1)=P_i=P_j$). Then, we get that
$$T(s)=\prod_{i=1}^nH_{P_i}(s)^{m_i-1}V(s),$$ where   $V(s)\in {\Bbb K}[s]$ and   $\gcd(H_{P_i},V)=1$ for
$i=1,\ldots,n$.

\para

Note that if  $V(s_0)=0$, then   $T(s_0)=0$ and thus,  $R_{12}(s_0)=R_{13}(s_0)=R_{23}(s_0)=0$. From $R_{13}(s_0)=\Res_t(G^*_1(s,t),G^*_3(s,t))(s_0)=0$ and using the properties of the resultant, we deduce that one of the following two statements hold:
\begin{enumerate}
\item There exists $s_1\in {\Bbb K}$ such that  $G^*_1(s_0,s_1)=G^*_3(s_0,s_1)=0$. Thus,  $H_P(s_1)=0$, where $P=\mathcal{P}(s_0)$, which is impossible since $\gcd(V,H_P)=1$.
\item It holds that  $\gcd(\lc_t(G^*_1),\,\lc_t(G_3^*))(s_0)=0$. However, this is also a contraction since we would have that
$$\lc_t(G^*_1)(s_0)=\lc_t(G_1)(s_0)=p_1(s_0)c_d-p(s_0)a_d=0\Rightarrow \frac{p_1(s_0)}{p(s_0)}=\frac{a_d}{c_d}$$
and
$$\lc_t(G^*_3)(s_0)=\lc_t(G_3)(s_0)=p_1(s_0)b_d-p_2(s_0)a_d=0\Rightarrow
\frac{p_1(s_0)}{p_2(s_0)}=\frac{a_d}{b_d}.$$ From both equalities, we deduce that
$$\frac{p_2(s_0)}{p(s_0)}=\frac{b_d}{c_d},$$ and then  $\mathcal{P}(s_0)=P_L$. This would imply that $P_L$ can be reached by the parametrization ${\cal P}(t)$ but this implies that it is a singularity (see Proposition 3.4 in \cite{MyB-2017(b)}), which contradicts the assumption of the theorem.\end{enumerate}

\para

  Thus, we have that  $V\in {\Bbb K}$ and, hence, we conclude that, up to constants in ${\Bbb K}\setminus\{0\}$,
$$T(s)=\prod_{i=1}^nH_{P_i}(s)^{m_i-1}.$$
\hfill $\Box$

\para

\end{document}